\pdfoutput=1
\RequirePackage{ifpdf}
\ifpdf 
\documentclass[pdftex]{sigma}
\else
\documentclass{sigma}
\fi

\numberwithin{equation}{section}

\newtheorem{thm}{Theorem}[section]
\newtheorem{cor}[thm]{Corollary}
\newtheorem{lem}[thm]{Lemma}
\newtheorem{prop}[thm]{Proposition}

\theoremstyle{definition}
\newtheorem*{rem}{Remark}
\newtheorem*{Example}{Example}

\begin{document}


\renewcommand{\thefootnote}{$\star$}

\renewcommand{\PaperNumber}{048}

\FirstPageHeading

\ArticleName{Combinatorial Formulae for Nested Bethe Vectors\footnote{This paper is a~contribution to the Special Issue in honor of Anatol Kirillov and Tetsuji Miwa.
The full collection is available at
\href{http://www.emis.de/journals/SIGMA/InfiniteAnalysis2013.html}{http://www.emis.de/journals/SIGMA/InfiniteAnalysis2013.html}}}

\ShortArticleName{Combinatorial Formulae for Nested Bethe Vectors}

\AuthorNameForHeading{V.~Tarasov  and  A.~Varchenko}

\Author{Vitaly TARASOV~$^{\dag\ddag}$ and  Alexander VARCHENKO~$^\S$}

\Address{$^\dag$~Department of Mathematical Sciences,
Indiana University~-- Purdue University Indianapolis,\\
\hphantom{$^\dag$}~402 North Blackford St, Indianapolis, IN 46202-3216, USA}
\EmailD{\href{mailto:vt@math.iupui.edu}{vt@math.iupui.edu}}

\Address{$^\ddag$~St.~Petersburg Branch of Steklov Mathematical Institute,\\
\hphantom{$^\ddag$}~Fontanka 27, St.~Petersburg, 191023, Russia}
\EmailD{\href{mailto:vt@pdmi.ras.ru}{vt@pdmi.ras.ru}}

\Address{$^\S$~Department of Mathematics, University of North Carolina
at Chapel Hill,\\
\hphantom{$^\S$}~Chapel Hill, NC 27599-3250, USA}
\EmailD{\href{mailto:anv@email.unc.edu}{anv@email.unc.edu}}

\ArticleDates{Received March 21, 2013, in f\/inal form June 27, 2013; Published online July 19, 2013}

\Abstract{We give combinatorial formulae for vector-valued weight functions (of\/f-shell nested Bethe vectors)
for tensor products of irreducible evaluation modules over the Yangian $Y({\mathfrak{gl}}_N)$ and the quantum af\/f\/ine algebra~$U_q(\widetilde{{\mathfrak{gl}}_N})$.}

\Classification{82B23; 17B80; 17B37; 81R50}

\renewcommand{\thefootnote}{\arabic{footnote}}
\setcounter{footnote}{0}

\pdfbookmark[1]{Introduction}{intr}
\section*{Introduction}
\label{intro}

In this paper we give combinatorial formulae for vector-valued weight functions for tensor \mbox{products}
of irreducible evaluation modules over the Yangian $Y({\mathfrak{gl}}_N)$ and the quantum af\/f\/ine algebra~$U_q(\widetilde{{\mathfrak{gl}}_N})$. Those functions are
also known as (of\/f-shell) nested Bethe vectors. They play an important role in
the theory of quantum integrable models and representation theory of Lie algebras and
quantum groups.

The nested algebraic Bethe ansatz was developed as a tool to f\/ind eigenvectors and eigenvalues
of transfer matrices of lattice integrable models associated with higher rank
Lie algebras, see~\cite{KR}. Similar to the regular Bethe ansatz, which is used in
the rank one case, eigenvectors are obtained as values of a certain rational function (nested
Bethe vector) on solutions of some system of algebraic equations (Bethe ansatz equations). Later, the nested
Bethe vectors (also called vector-valued weight functions) were used to construct Jackson integral representations for
solutions of the quantized (dif\/ference) Knizhnik--Zamolodchikov (qKZ) equations~\cite{TV1}. The results
of~\cite{KR} has been extended to higher transfer matrices in \cite{MTV1}.

In the rank one case combinatorial formulae for vector-valued weight function are important
in various areas from computation of correlation functions in integrable models,
see~\cite{KBI}, to eva\-luation of some multidimensional generalizations of
the Vandermonde determinant \cite{TV2}. In the ${\mathfrak{gl}}_N$ case considered in
this paper, combinatorial formulae, in particular, clarify analytic properties
of the vector-valued weight function, which is important for constructing hypergeometric solutions of
the qKZ equations associated with ${\mathfrak{gl}}_N$.

Combinatorial formulae for the vector-valued weight functions associated with the dif\/ferential Knizhnik--Zamolodchikov equations
were developed in \cite{M,SV1,SV2,RSV, FRV}.

The results of this paper were obtained while the authors were visiting
the Max-Planck-Institut f\"ur Mathematik in Bonn in 1998.
The results of the paper were used in \cite{MTT,KPT,TV3}.

The paper has appeared in the arXiv in 2007, but still looks topical.
It is published with no intention to give any review of the subject or ref\/lect
the state of the art. Let us only mention a~few papers making progress in
particularly close problems \cite{KP1,KP2,FKPR,OPS, BPRS1,BPRS2} and those
exploring recently the results of the paper \cite{RTV,TV4, MTV2}.

The paper is organized as follows. In Sections \ref{:Rwtf}--\ref{:SBBvox}
we consider in detail the Yangian case. In the traditional terminology this
case is called \textit{rational}. In Section~\ref{:Twtf} we formulate the results
for the quantum af\/f\/ine algebra case, also called \textit{trigonometric}. The proofs in that case are
very similar to the Yangian case.

\section{Basic notation}
We will be using the standard superscript notation for embeddings of tensor
factors into tensor products. If ${\mathcal A}_1,\dots,{\mathcal A}_k$ are unital associative
algebras, and $a\in {\mathcal A}_i$, then
\[
a^{(i)}= 1^{\otimes(i-1)}\otimes a \otimes 1^{\otimes(k-i)} \in {\mathcal A}_1\otimes \cdots \otimes {\mathcal A}_k .
\]
If $a\in {\mathcal A}_i$ and $b\in {\mathcal A}_j$, then $(a\otimes b)^{(ij)} =a^{(i)}b^{(j)}$,
etc.
\begin{Example}
Let $k=2$. Let ${\mathcal A}_1$, ${\mathcal A}_2$ be two copies of the same algebra ${\mathcal A}$.
Then for any $a,b\in {\mathcal A}$ we have $a^{(1)}=a\otimes1$, $b^{(2)}=1\otimes b$,
$(a\otimes b)^{(12)}= a\otimes b$ and $(a\otimes b)^{(21)}= b\otimes a$.
\end{Example}

Fix a positive integer $N$. All over the paper we identify elements of
$\operatorname{End}\big({\mathbb C}^N\big)$ with ${N\times N}$ matrices using the standard basis of ${\mathbb C}^N$.

We will use the rational and trigonometric $R$-matrices. The \textit{rational $R$-matrix} is
\begin{gather}
R(u) = u + \sum_{a,b=1}^N E_{ab}\otimes E_{ba} ,
\label{R}
\end{gather}
where $E_{ab}\in\operatorname{End}\big({\mathbb C}^N\big)$ is a matrix with the only nonzero entry equal
to $1$ at the intersection of the  $a$-th row and $b$-th column.
The $R$-matrix satisf\/ies the inversion relation
\[
R(u) R^{(21)}(-u)=1-u^2
\]
and the Yang--Baxter equation
\begin{gather}
R^{(12)}(u-v) R^{(13)}(u) R^{(23)}(v) =
R^{(23)}(v) R^{(13)}(u) R^{(12)}(u-v) .
\label{YB}
\end{gather}

Fix a complex number $q$ not equal to $\pm1$. The \textit{trigonometric $R$-matrix}
\begin{gather}
R_q(u) = \big(u q-q^{-1}\big) \sum_{a=1}^N E_{aa}\otimes E_{aa} +
(u-1) \sum_{1\leqslant a<b\leqslant N} (E_{aa}\otimes E_{bb}+E_{bb}\otimes E_{aa})
\nonumber\\
\hphantom{R_q(u) =}{}
+ \big(q-q^{-1}\big) \sum_{1\leqslant a<b\leqslant N} (u E_{ab}\otimes E_{ba}+E_{ba}\otimes E_{ab})\label{Rq}
\end{gather}
satisf\/ies the inversion relation
\[
R_q(u) R_q^{(21)}\big(u^{-1}\big)= \big(u q-q^{-1}\big) \big(u^{-1}q-q^{-1}\big)
\]
and the Yang--Baxter equation
\[
R_q^{(12)}(u/v) R_q^{(13)}(u) R_q^{(23)}(v)=
R_q^{(23)}(v) R_q^{(13)}(u) R_q^{(12)}(u/v).
\]

Let $e_{ab}$, $a,b=1,\dots,N$, be the standard generators of the Lie algebra
${\mathfrak{gl}}_N$:
\[
[e_{ab},e_{cd}]=\delta_{bc} e_{ad}-\delta_{ad} e_{cb}.
\]
Let ${\mathfrak{h}} =\bigoplus_{a=1}^N{\mathbb C} e_{aa} $ be the Cartan subalgebra.
For any $\Lambda\in{\mathfrak{h}}^*$ we set $\Lambda^{a} =\langle\Lambda ,e_{aa}\rangle$,
and identify~${\mathfrak{h}}^*$ with ${\mathbb C}^N $ by taking $\Lambda$ to $\big(\Lambda^{1},\dots,\Lambda^{N}\big)$.
We use the Gauss decomposition ${\mathfrak{gl}}_N = {\mathfrak{h}}\oplus{\mathfrak{n}}_{+} \oplus{\mathfrak{n}}_{-}$
where ${\mathfrak{n}}_{+}=\bigoplus_{a<b}{\mathbb C} e_{ab}$ and ${\mathfrak{n}}_{-}=\bigoplus_{a<b}{\mathbb C} e_{ba}$.
A vector $v$ in a ${\mathfrak{gl}}_N$-module is called a~\textit{singular vector} if ${\mathfrak{n}}_{+} v=0$.
The space ${\mathbb C}^N $ is considered as a~${\mathfrak{gl}}_N$-module with the natural action,
$e_{ab} \mapsto E_{ab}$. This module is called the \textit{vector representation}.

\section{Rational weight functions}\label{:Rwtf}

The Yangian $Y({\mathfrak{gl}}_N)$ is a unital associative algebra with generators
$T_{ab}^{\{s\}}$, $a,b=1,\dots,N$, and $s=1,2,\ldots$.
Organize them into generating series:
\begin{gather}
T_{ab}(u)=\delta_{ab}+\sum_{s=1}^\infty T_{ab}^{\{s\}} u^{-s} ,\qquad a,b=1,\dots,N.
\label{Tab}
\end{gather}
The def\/ining relations in $Y({\mathfrak{gl}}_N)$ have the form
\begin{gather}
(u-v)\bigl[T_{ab}(u),T_{cd}(v)\bigr]=
T_{cb}(v)T_{ad}(u)-T_{cb}(u) T_{ad}(v),
\label{Tabcd}
\end{gather}
for all $a,b,c,d=1,\dots,N$.

Combine series \eqref{Tab} together into a series
${T(u)=\sum\limits_{a,b=1}^{N}E_{ab}\otimes T_{ab}(u)}$ with coef\/f\/icients in
$\operatorname{End}\big({\mathbb C}^N\big)\otimes Y({\mathfrak{gl}}_N)$. Relations~\eqref{Tabcd} amount to the following
equality for series with coef\/f\/icients in  $\operatorname{End}\big({\mathbb C}^N\big)\otimes\operatorname{End}\big({\mathbb C}^N\big)\otimes Y({\mathfrak{gl}}_N)$:
\begin{gather}
R^{(12)}(u-v) T^{(13)}(u) T^{(23)}(v)=
T^{(23)}(v) T^{(13)}(u) R^{(12)}(u-v).
\label{RTT}
\end{gather}

The Yangian $Y({\mathfrak{gl}}_N)$ is a Hopf algebra. In terms of generating series~\eqref{Tab},
the coproduct $\Delta:Y({\mathfrak{gl}}_N)\to Y({\mathfrak{gl}}_N)\otimes Y({\mathfrak{gl}}_N)$ reads as follows:
\begin{gather}
\Delta\bigl(T_{ab}(u)\bigr)=\sum_{c=1}^N T_{cb}(u)\otimes T_{ac}(u),
\qquad a,b=1,\dots,N.
\label{Dl}
\end{gather}

There is a one-parameter family of automorphisms ${\rho_x:Y({\mathfrak{gl}}_N)\to Y({\mathfrak{gl}}_N)}$ def\/ined in
terms of the series $T(u)$ by the rule ${\rho_x\bigl(T(u)\bigr)= T(u-x)}$;
in the right side, $(u-x)^{-1}$ has to be expanded as a power series in $u^{-1}$.

The Yangian $Y({\mathfrak{gl}}_N)$ contains the universal enveloping algebra $U({\mathfrak{gl}}_N)$ as a Hopf
subalgebra. The embedding is given by $e_{ab}\mapsto T_{ba}^{\{1\}}$ for all
$a,b=1,\dots,N$. We identify $U({\mathfrak{gl}}_N)$ with its image in $Y({\mathfrak{gl}}_N)$ under this embedding.
It is clear from relations~\eqref{Tabcd} that for any $a,b=1,\dots,N$,
\[
\bigl[E_{ab}\otimes1+1\otimes e_{ab},T(u)\bigr]=0.
\]

The \textit{evaluation homomorphism} ${\epsilon: Y({\mathfrak{gl}}_N)\to U({\mathfrak{gl}}_N)}$ is given by the rule
${\epsilon: T_{ab}^{\{1\}} \mapsto e_{ba}}$ for any $a,b=1,\dots,N$, and
$\epsilon: T_{ab}^{\{s\}} \mapsto 0$ for any $s>1$ and all $a$, $b$.
Both the automorphisms $\rho_x$ and the homomorphism $\epsilon$ restricted to the subalgebra
$U({\mathfrak{gl}}_N)$ are the identity maps.

For a ${\mathfrak{gl}}_N$-module $V$ denote by $V(x)$ the $Y({\mathfrak{gl}}_N)$-module induced from $V$ by the homomorphism
$\epsilon\circ \rho_x$. The module $V(x)$ is called an \textit{evaluation module} over $Y({\mathfrak{gl}}_N)$.

A vector $v$ in a $Y({\mathfrak{gl}}_N)$-module is called \textit{singular} with respect to the action of $Y({\mathfrak{gl}}_N)$
if $T_{ba}(u) v = 0$ for all $1\leqslant a<b\leqslant N$. A singular vector $v$ that is an eigenvector
for the action of $T_{11}(u),\dots, T_{N N}(u)$ is called a \textit{weight singular
vector}; the respective eigenvalues are denoted by
$\langle T_{11}(u) v\rangle,\dots, \langle T_{N N}(u) v\rangle$.

\begin{Example}
Let $V$ be a ${\mathfrak{gl}}_N$-module and let ${v\in V}$ be a singular vector of weight
$\big(\Lambda^{1},\dots,\Lambda^{N}\big)$. Then~$v$ is a weight singular vector with respect to the action of $Y({\mathfrak{gl}}_N)$
in the evaluation mo\-du\-le~$V(x)$ and $\langle T_{aa}(u) v\rangle = 1+\Lambda^{a}(u-x)^{-1}$,
$a=1,\dots,N$.
\end{Example}

If $v_1$, $v_2$ are weight singular vectors with respect to the action of $Y({\mathfrak{gl}}_N)$ in $Y({\mathfrak{gl}}_N)$-modules
$V_1$, $V_2$, then the vector ${v_1 \otimes v_2}$ is a weight singular vector
with respect to the action of $Y({\mathfrak{gl}}_N)$ in the tensor product $V_1\otimes V_2$, and
$\langle T_{aa}(u) v_1 \otimes v_2\rangle = \langle T_{aa}(u) v_1\rangle
\langle T_{aa}(u) v_2\rangle$ for all $a=1,\dots,N$.

We will use two embeddings of the algebra $ Y({\mathfrak{gl}}_{N-1})$ into $Y({\mathfrak{gl}}_N)$,
called $\phi$ and $\psi$:
\begin{gather}
\phi \bigl(T_{ab}^{\langle N-1\rangle}(u)\bigr) = T_{ab}^{\langle N\rangle}(u) ,\qquad
\psi \bigl(T_{ab}^{\langle N-1\rangle }(u)\bigr) = T_{a+1, b+1}^{\langle N\rangle} (u),
\label{embed}
\end{gather}
$a,b=1,\dots,N-1$. Here $T_{ab}^{\langle N-1\rangle}(u)$ and $T_{ab}^{\langle N\rangle}(u)$ are series~\eqref{Tab}
for the algebras $ Y({\mathfrak{gl}}_{N-1})$ and $Y({\mathfrak{gl}}_N)$, respectively.

Let $\xi=(\xi^1,\dots,\xi^{N-1})$ be a collection of nonnegative integers.
Set ${\xi^{<a}=\xi^1+\cdots +\xi^{a-1}}$, $a=1,\dots,N$, and
${|\xi|=\xi^1+\cdots +\xi^{N-1}=\xi^{<N}}$.
Consider a series in $|\xi|$ variables $t^1_1,\dots, t^1_{\xi^1}$, $\dots, t^{N-1}_1,\dots,t^{N-1}_{\xi^{N-1}}$ with coef\/f\/icients in $Y({\mathfrak{gl}}_N)$:
\begin{gather}
{\widehat{\mathbb B}}_{ \xi}\big(t^1_1,\dots, t^{N-1}_{\xi^{N-1}}\big) = \big({\operatorname{tr}^{\otimes|\xi|}} \otimes {\rm id}\big) \biggl(
T^{(1, |\xi|+1)}\big(t^1_1\big)\cdots T^{(|\xi|,|\xi|+1)}\big(t^{N-1}_{\xi^{N-1}}\big)
\nonumber\\
\phantom{{\widehat{\mathbb B}}_{ \xi}(t^1_1,\dots, t^{N-1}_{\xi^{N-1}}) =}{}
\times \mathop{\overrightarrow\prod}\limits_{(a,i)<(b,j)}
R^{(\xi^{<b}+j,\xi^{<a}+i)}\big(t^b_j -t^{a}_i\big)
E_{21}^{\otimes\xi^1} \otimes\cdots \otimes  E_{N,N-1}^{\otimes\xi^{N-1}}\otimes 1 \biggr).\label{BBh}
\end{gather}
Here ${\operatorname{tr}:\operatorname{End}\big({\mathbb C}^N\big)\to{\mathbb C}}$ is the standard trace map, the pairs in the product
are ordered lexi\-co\-graphically, $(a,i)<(b,j)$ if $a<b$, or $a=b$
and $i<j$; the product is taken over all two-element subsets
of the set $\big\{(c,k) \vert c=1,\dots, N-1 ,\  k=1,\dots,\xi^c \big\}$;
in the product,
the factor $R^{(\xi^{<b}+j,\xi^{<a}+i)}\big(t^b_j-t^{a}_i\big)$ is
to the left of $R^{(\xi^{<d}+ l,\xi^{<c}+k)}\big(t^d_l-t^c_k\big)$
if $(a,i)<(c,k)$, or $(a,i)=(c,k)$ and $(b,j)<(d,l)$.

\begin{Example}
Let $N=4$ and $\xi=(1,1,1)$. Then
\begin{gather*}
{\widehat{\mathbb B}}_{\xi}\big(t^1_1,t^2_1,t^3_1\big) = \big({\operatorname{tr}^{\otimes 3}} \otimes {\rm id}\big)\Bigl(
T^{(14)}\big(t^1_1\big) T^{(24)}\big(t^2_1\big) T^{(34)}\big(t^3_1\big)
\\
\hphantom{{\widehat{\mathbb B}}_{\xi}\big(t^1_1,t^2_1,t^3_1\big) =}{}
 \times R^{(32)}\big(t^3_1-t^2_1\big) R^{(31)}\big(t^3_1 -t^1_1\big)
R^{(21)}\big(t^2_1 -t^1_1\big) E_{21}\otimes E_{32}\otimes E_{43}\otimes 1\Bigr).
\end{gather*}
\end{Example}

\begin{rem}
The series ${\widehat{\mathbb B}}_{\xi}\big(t^1_1,\dots, t^{N-1}_{\xi^{N-1}}\big)$ belongs to  ${Y({\mathfrak{gl}}_N)\bigl[t^1_1,\dots, t^{N-1}_{\xi^{N-1}}\bigr]\bigl[\bigl[\big(t^1_1\big)^{-1}\!,\dots,\big(t^{N-1}_{\xi^{N-1}}\big)^{-1}\bigr]\bigr]}$.
\end{rem}

\begin{rem}
Relations \eqref{RTT} imply that
\begin{gather}
 T^{(1,|\xi|+1)}\big(t^1_1\big)\cdots T^{(|\xi|,|\xi|+1)}\big(t^{N-1}_{\xi^{N-1}}\big)
\mathop{\overrightarrow\prod}\limits_{(a,i)<(b,j)}
R^{(\xi^{<b}+ j,\xi^{<a}+ i)}\big(t^b_j-t^{a}_i\big)
\nonumber
\\
\qquad {} = \mathop{\overrightarrow\prod}\limits_{(a,i)<(b,j)}
R^{(\xi^{<b}+j,\xi^{<a}+i)}\big(t^b_j-t^{a}_i\big)
T^{(|\xi|,|\xi|+1)}\big(t^{N-1}_{\xi^{N-1}}\big)\cdots T^{(1,|\xi|+1)}\big(t^1_1\big).\label{TTR}
\end{gather}
For instance,
\begin{gather*}
T^{(14)}\big(t^1_1\big) T^{(24)}\big(t^2_1\big) T^{(34)}\big(t^3_1\big)
R^{(32)}\big(t^3_1-t^2_1\big) R^{(31)}\big(t^3_1-t^1_1\big) R^{(21)}\big(t^2_1-t^1_1\big)
\\
\qquad {}= R^{(32)}\big(t^3_1-t^2_1\big) R^{(31)}\big(t^3_1-t^1_1\big)
R^{(21)}\big(t^2_1-t^1_1\big) T^{(34)}\big(t^3_1\big) T^{(24)}\big(t^2_1\big) T^{(14)}\big(t^1_1\big).
\end{gather*}
\end{rem}

\begin{rem}
Using the Yang--Baxter equation \eqref{YB} one can rearrange the factors in the product of $R$-matrices
in formulae~\eqref{BBh}, \eqref{TTR}. For instance,
\[
\mathop{\overrightarrow\prod}\limits_{(a,i)<(b,j) }
R^{(\xi^{<b} + j, \xi^{<a} + i)}\big(t^b_j -t^{a}_i\big) =
\mathop{\overleftarrow\prod}\limits_{(a,i)<(b,j)}
R^{(\xi^{<b} + j, \xi^{<a} + i)}\big(t^b_j -t^{ a}_i\big),
\]
where in the right side the factor
$R^{(\xi^{<b} + j, \xi^{<a} + i)}\big(t^b_j -t^{ a}_i\big)$ is to the right of
$R^{(\xi^{<d} + l,\xi^{<c}+ k)}\big(t^d_l -t^c_k\big)$ if $(a,i)<(c ,k)$,
or $(a,i)=(c,k)$ and $(b,j)<(d,l)$. In particular, for any $a=1,\dots,N-1$, and
any $i=1,\dots,\xi^a -1$, there are rearrangements of factors such that
$R^{(\xi^{<a} + i +1, \xi^{<a}+i)}\big(t^{a}_{i+1}-t^{a}_i\big)$
is the leftmost or the rightmost factor of the product.
\end{rem}

\begin{Example}
Let $N=4$ and $\xi=(2,1,1)$. Then
\begin{gather*}
  R^{(43)}\big(t^3_1 -t^2_1\big) R^{(42)}\big(t^3_1 -t^1_2\big) R^{(41)}\big(t^3_1 -t^1_1\big)
R^{(32)}\big(t^2_1 -t^1_2\big) R^{(31)}\big(t^2_1 -t^1_1\big) R^{(21)}\big(t^1_2 -t^1_1\big)
\\
\qquad\! {}=
R^{(32)}\big(t^2_1 -t^1_2\big) R^{(31)}\big(t^2_1 -t^1_1\big) R^{(21)}\big(t^1_2 -t^1_1\big)
R^{(41)}\big(t^3_1 -t^1_1\big) R^{(42)}\big(t^3_1-t^1_2\big) R^{(43)}\big(t^3_1 -t^2_1\big).
\end{gather*}
\end{Example}

Further on, we will abbreviate, $t=\big(t^1_1,\dots, t^{N-1}_{\xi^{N-1}}\big)$. Set
\begin{gather}
{\mathbb B}_{\xi}(t) = {\widehat{\mathbb B}}_{\xi}(t)
\prod_{a=1}^{N-1} \prod_{1\leqslant i<j\leqslant\xi^a}\frac1{t^{a}_j -t^{a}_i+1}
\prod_{1\leqslant a<b<N}\prod_{i=1}^{\xi^a}\prod_{j=1}^{\xi^b}
\frac1{t^b_j -t^{a}_i} ,
\label{BB}
\end{gather}
cf.~\eqref{BBh}. To indicate the dependence on $N$, if necessary,
we will write ${\mathbb B}^{\langle N\rangle}_{\xi}(t)$.

\begin{Example}
Let $N=2$ and $\xi=(\xi^1)$.
Then  ${\mathbb B}^{\langle 2\rangle}_{\xi}(t) = T_{12}\big(t^1_1\big)\cdots T_{12}\big(t^1_{\xi^1}\big)$.
\end{Example}

\begin{Example}
Let $N=3$ and $\xi=(1,1)$. Then
\begin{gather*}
{\mathbb B}^{\langle 3\rangle}_{\xi}(t) = T_{12}\big(t^1_1\big) T_{23}\big(t^2_1\big)\! +
\frac1{t^2_1 -t^1_1} T_{13}\big(t^1_1\big) T_{22}\big(t^2_1\big)
= T_{23}\big(t^2_1\big) T_{12}\big(t^1_1\big)\! +
\frac1{t^2_1 -t^1_1} T_{13}\big(t^2_1\big) T_{22}\big(t^1_1\big).
\end{gather*}
\end{Example}

\begin{Example}
Let $N=4$ and $\xi=(1,1,1)$. Then
\begin{gather*}
{\mathbb B}^{\langle 4\rangle}_{\xi}(t) = T_{12}\big(t^1_1\big) T_{23}\big(t^2_1\big) T_{34}\big(t^3_1\big) \\
\hphantom{{\mathbb B}^{\langle 4\rangle}_{\xi}(t) =}{}
+\frac1{t^2_1-t^1_1} T_{13}\big(t^1_1\big) T_{22}\big(t^2_1\big) T_{34}\big(t^3_1\big)
 + \frac1{t^3_1 -t^2_1} T_{12}\big(t^1_1\big) T_{24}\big(t^2_1\big) T_{33}\big(t^3_1\big)
 \\
\hphantom{{\mathbb B}^{\langle 4\rangle}_{\xi}(t) =}{}
  +
\frac1{\big(t^2_1 -t^1_1\big) \big(t^3_1 -t^2_1\big)}
\bigl( T_{14}\big(t^1_1\big) T_{22}\big(t^2_1\big) T_{33}\big(t^3_1\big)
+
T_{13}\big(t^1_1\big) T_{24}\big(t^2_1\big) T_{32}\big(t^3_1\big)\bigr)
\\
\hphantom{{\mathbb B}^{\langle 4\rangle}_{\xi}(t) =}{}
+
\frac{\big(t^2_1 -t^1_1\big)\big(t^3_1 -t^2_1\big)+1}
{\big(t^2_1 -t^1_1\big)\big(t^3_1 -t^1_1\big)\big(t^3_1 -t^2_1\big)}
T_{14}\big(t^1_1\big) T_{23}\big(t^2_1\big) T_{32}\big(t^3_1\big)  .
\end{gather*}
\end{Example}

The direct product of the symmetric groups $S_{\xi^1}\times \cdots \times S_{\xi^{N-1}}$ acts on expressions in $|\xi|$
variables, permuting the variables with the same superscript:
\begin{gather}
\sigma^1 \times\cdots\times \sigma^{N-1}:\ f\big(t^1_1,\dots, t^{N-1}_{\xi^{N-1}}\big) \mapsto f\big(t^1_{\sigma^{1}_1},\dots, t^1_{\sigma^{1}_{\xi^1}}; \ldots ; t^{N-1}_{\sigma^{N-1}_1},\dots, t^{N-1}_{\sigma^{N-1}_{\xi^{N-1}}}\big) ,
\label{sixi}
\end{gather}
where  ${\sigma^a \in S_{\xi^a}}$,  $a=1,\dots,N-1$.

\pagebreak

\begin{lem}[\protect{\cite[Theorem~3.3.4]{TV1}}]
\label{BBS}
The expression ${\mathbb B}_{\xi}(t)$ is invariant under the action
of the group $S_{\xi^1}\times \cdots \times S_{\xi^{N-1}}$.
\end{lem}

\begin{proof}
Let $P =\sum_{a,b}E_{ab} \otimes E_{ba}$ be the f\/lip map, and ${{\check R}(u)=PR(u)}$.
For any $a=1,\dots,N-1$, we have
\begin{gather}
{\check R}(u) E_{a+1,a} \otimes E_{a+1,a} = (u+1) E_{a+1,a} \otimes E_{a+1,a} =
E_{a+1,a} \otimes E_{a+1,a} {\check R}(u) .
\label{REE}
\end{gather}

Set
\[
{\mathbb T}(t) =
T^{(1, |\xi|+1)}\big(t^1_1\big)\cdots T^{(|\xi|, |\xi|+1)}\big(t^{N-1}_{\xi^{N-1}}\big)
\mathop{\overrightarrow\prod}\limits_{(a,i)<(b,j)}
R^{(\xi^{<b} + j, \xi^{<a} + i }\big(t^b_j -t^{a}_i\big) .
\]
Let $\tilde t=\big({\tilde t}^{1}_{1},\dots,{\tilde t}^{N-1}_{\xi^{N-1}}\big) $ be obtained
from $t=\big(t^1_1,\dots, t^{N-1}_{\xi^{N-1}}\big)$ by the permutation of $t^{a}_i$ and $t^{a}_{i+1}$.
Set $j=i+\sum\limits_{b<a} \xi^{b}$. The Yang--Baxter equation \eqref{YB} and relations \eqref{RTT}
yield
\[
{\mathbb T}(t) {\check R}^{(j+1,j)}\big(t^{a}_i-t^{a}_{i+1}\big) =
{\check R}^{(j,j+1)}\big(t^{a}_{i+1}-t^{a}_i\big) {\mathbb T}(\tilde t).
\]
Hence,
\begin{gather*}
 {\widehat{\mathbb B}}_{\xi}(t)=\big({\operatorname{tr}^{\otimes|\xi|}} \otimes {\rm id}\big) \Big({\mathbb T}(t)
E_{21}^{\otimes\xi^1}\otimes\cdots\otimes E_{N,N-1}^{\otimes\xi^{N-1}}\otimes 1 \Big)
\\
 = \big({\operatorname{tr}^{\otimes|\xi|}} \otimes {\rm id}\big) \Big(
{\check R}^{(j,j+1)}\big(t^{a}_{i+1}-t^{a}_i\big){\mathbb T}(\tilde t)
\bigl({\check R}^{(j+1,j)}\big(t^{a}_i-t^{a}_{i+1}\big)\bigr)^{-1}
E_{21}^{\otimes\xi^1} \otimes\cdots\otimes E_{N,N-1}^{\otimes\xi^{N-1}}\otimes 1 \Big)
\\
 = \frac{t^{a}_{i+1}-t^{a}_i+1}{t^{a}_i-t^{a}_{i+1}+1}
\big({\operatorname{tr}^{\otimes|\xi|}} \otimes {\rm id}\big) \Big({\mathbb T}(\tilde t)
E_{21}^{\otimes\xi^1}\otimes\cdots\otimes E_{N,N-1}^{\otimes\xi^{N-1}}\otimes 1\Big)
=\frac{t^{a}_{i+1}-t^{a}_i+1}{t^{a}_i-t^{a}_{i+1}+1}
 {\widehat{\mathbb B}}_{\xi}(\tilde t),
\end{gather*}
by formula \eqref{REE} and the cyclic property of the trace.
Therefore, ${\mathbb B}_{\xi}(\tilde t)={\mathbb B}_{\xi}(t)$, see~\eqref{BB}.
\end{proof}

If $v$ is a weight singular vector with respect to the action of $Y({\mathfrak{gl}}_N)$,
we call the expression ${\mathbb B}_{\xi}(t)v$ the (\textit{rational})
\textit{vector-valued weight function} of weight $\big(\xi^1,{\xi^2-\xi^1},\dots,
{\xi^{N-1}-\xi^{N-2}}$, $- \xi^{N-1}\big)$ associated with~$v$.

Weight functions associated with ${\mathfrak{gl}}_N$ weight singular vectors in evaluation $Y({\mathfrak{gl}}_N)$-modules
(in particular, highest weight vectors of highest weight ${\mathfrak{gl}}_N$-modules) can be calculated explicitly
by means of the following Theorems~\ref{BBN} and \ref{BB1}. The theorems express
weight functions for $Y({\mathfrak{gl}}_N)$ in terms of weight functions for $ Y({\mathfrak{gl}}_{N-1})$. Applying the theorems
several times one can get $2^{N-2}$ combinatorial expressions for the same
weight function, the expressions being labeled by subsets of $\{1,\dots, N-2\}$.
The expressions corresponding to the empty set and the whole set are given
in Corollaries~\ref{BCN} and~\ref{BC1}.

Let $v_1,\dots,v_n$ be weight singular vectors with respect to the action of $Y({\mathfrak{gl}}_N)$.
Corollary~\ref{BBvoxn} expresses the weight function ${{\mathbb B}_{\xi}(t) (v_1\otimes\cdots\otimes v_n)}$
as a sum of the tensor products
${{\mathbb B}_{\zeta_1} (t_1) v_1\otimes\cdots\otimes{\mathbb B}_{\zeta_n} (t_n) v_n}$ with
${\zeta_1 +\cdots+\zeta_n = \xi}$, and $t_1,\dots, t_n$ being a partition of
the collection $t$ of $|\xi|$ variables into collections of
$| \zeta_1|,\dots, |\zeta_n|$ variables. This yields combinatorial formulae for
weight functions associated with tensor products of highest weight vectors of highest weight evaluation modules.

\begin{rem}
It is shown in \cite{KR} that for a weight singular vector $v$ in a tensor
product of evaluation  $Y({\mathfrak{gl}}_N)$-modules, the values of the weight function ${\mathbb B}_{\xi}(t) v$ at
solutions of a certain system of algebraic equations (Bethe ansatz equations) are eigenvectors of the
transfer matrix of the corresponding lattice integrable model. This result
is extended in~\cite{MTV1} to the case of higher transfer matrices.
\end{rem}
\begin{rem}
The weight functions ${\mathbb B}_{\xi}(t) v$ are used in~\cite{TV1} to construct
Jackson integral representations for solutions of the qKZ equations.
\end{rem}

\begin{rem}
The expression for a vector-valued weight function used here may dif\/fer from
the expressions for the corresponding objects used in other papers,
see~\cite{KR,TV1}. The discrepancy is not essential and may occur
due to the choice of coproduct for the Yangian $Y({\mathfrak{gl}}_N)$ as well as the choice of
normalization.
\end{rem}

\section{Combinatorial formulae for rational weight functions}
For a nonnegative integer $k$ introduce a function $W_{k}(t_1,\dots,t_k)$:
\[
W_{ k}(t_1,\dots,t_k) =
\prod_{1\leqslant i<j\leqslant k} \frac{t_i -t_j -1}{t_i -t_j} .
\]
For an expression $f\big(t^1_1,\dots, t^{N-1}_{\xi^{N-1}}\big)$, set
\begin{gather}
\operatorname{Sym}^\xi_{ t}f\big(t^1_1,\dots, t^{N-1}_{\xi^{N-1}}\big) =
\sum_{\sigma^1,\dots, \sigma^{N-1}} f\big(t^1_{\sigma^{1}_1},\dots, t^1_{\sigma^{1}_{\xi^1}}; \ldots ; t^{N-1}_{\sigma^{N-1}_1},\dots, t^{N-1}_{\sigma^{N-1}_{\xi^{N-1}}}\big) ,
\label{Sym}
\end{gather}
where  ${\sigma^a \in S_{\xi^a}}$,  $a=1,\dots,N-1$,  and
\begin{gather}
\overline{\operatorname{Sym}}^\xi_{t}f(t) =\operatorname{Sym}^\xi_{t}\left(f(t)
\prod_{a=1}^{N-1}W_{\xi^a} \big(t^{a}_1,\dots, t^{a}_{\xi^a}\big)\right)  .
\label{Symb}
\end{gather}

Let ${\eta^1 \leqslant\cdots \leqslant\eta^{N-1}} $ be nonnegative integers.
Def\/ine a function $X_\eta\big(t^1_1,\dots, t^1_{\eta^1};\dots; t^{N-1}_1,\dots$, $t^{N-1}_{\eta^{N-1}}\big)$,
\begin{gather}
X_\eta(t)=\prod_{a=1}^{N-2}\left[
\prod_{j=1}^{\eta^a} \frac1{t^{a+1}_j-t^{a}_j}\prod_{i=1}^{j-1}
\frac{t^{a+1}_i-t^{a}_j+1}{t^{a+1}_i-t^{a}_j}\right] .
\label{X}
\end{gather}
The function $X_\eta(t)$ does not actually depend on
the variables $t^{N-1}_{\eta^{N-2}+1},\dots, t^{N-1}_{\eta^{N-1}}$.

For nonnegative integers ${\eta^{1}\geqslant\cdots\geqslant \eta^{N-1}}$
def\/ine a function $Y_{\eta}\big(t^1_1,\dots, t^1_{\eta^1};\dots; t^{N-1}_1,\dots$, \linebreak $t^{N-1}_{\eta^{N-1}}\big)$,
\begin{gather}
Y_{\eta}(t)=\prod_{a=2}^{N-1}\left[ \prod_{j=1}^{\eta^a}
 \frac1{t^{a}_j-t^{a-1}_{j+\eta^{a-1}-\eta^a}}
\prod_{i=1}^{j-1}\frac{t^{a}_i-t^{a-1}_{j+\eta^{a-1}-\eta^a}+1}
{t^{a}_i-t^{a-1}_{j+\eta^{a-1}-\eta^a}}\right] .
\label{Y}
\end{gather}
The function $Y_{\eta}(t)$ does not actually depend on the variables
$t^1_1,\dots, t^1_{\eta^1-\eta^2}$.

For any $\xi,\eta\in{\mathbb Z}_{\geqslant 0}^{N-1}$,
def\/ine a function $Z_{\xi,\eta}\big(t^1_1,\dots, t^{N-1}_{\xi^{N-1}};s^1_1,\dots,s^{N-1}_{\eta^{N-1}}\big)$,
\begin{gather}
Z_{\xi,\eta}(t;s)=\prod_{a=1}^{N-2}\prod_{i=1}^{\xi^{a+1}}
\prod_{j=1}^{\eta^a}\frac{t^{a+1}_i-s^a_j+1}{t^{a+1}_i-s^a_j}.
\label{Z}
\end{gather}
The function $Z_{\xi,\eta}(t;s)$ does not actually depend
on the variables $t^1_1,\dots, t^1_{\xi^1}$ and $s^{N-1}_1,\dots, s^{N-1}_{\eta^{N-1}}$.

If $\xi,\eta,\zeta\in{\mathbb Z}_{\geqslant 0}^{N-1}$ are such that $\xi-\zeta\in{\mathbb Z}_{\geqslant 0}^{N-1}$ and
$\zeta-\eta\in{\mathbb Z}_{\geqslant 0}^{N-1}$, and $t=(t^1_1,\dots, t^1_{\xi^1},\dots, t^{N-1}_1$, $\dots,t^{N-1}_{\xi^{N-1}})$, then we set{\samepage
\begin{gather}
t_{[\eta]}= \big(t^1_1,\dots, t^1_{\eta^1};\dots; t^{N-1}_1,\dots, t^{N-1}_{\eta^{N-1}}\big),
\nonumber
\\
t_{(\eta,\zeta]}=\big(t^1_{\eta^1+1},\dots, t^1_{\zeta^1};\dots;
t^{N-1}_{\eta^{N-1}+1},\dots, t^{N-1}_{\zeta^{N-1}}\big) . \label{teta}
\end{gather}
Notice that $t_{[\eta]}=t_{(0,\eta]}$.}

For any $\xi=\big(\xi^1,\dots,\xi^{N-1}\big)$ set $\dot\xi=\big(\xi^1,\dots,\xi^{N-2}\big)$ and
$\ddot\xi=\big(\xi^2,\dots,\xi^{N-1}\big)$. If $t=\big(t^1_1,\dots, t^1_{\xi^1}$, $\dots, t^{N-1}_1,\dots,t^{N-1}_{\xi^{N-1}}\big)$, then we set
\begin{gather}
\label{dot}
\dot t = \big(t^1_1,\dots, t^1_{\xi^1};\dots;t^{N-2}_1,\dots, t^{N-2}_{\xi^{N-2}}\big),
\qquad
\ddot t= \big(t^2_1,\dots, t^2_{\xi^2};
\dots;t^{N-1}_1,\dots, t^{N-1}_{\xi^{N-1}}\big).
\end{gather}

\begin{thm}\label{BBN}
Let $V$ be a ${\mathfrak{gl}}_N$-module and ${v\in V}$ a singular vector of weight $\big(\Lambda^{1},\dots,\Lambda^{N}\big)$.
Let $\xi^1,\dots,\xi^{N-1}$ be nonnegative integers and $t=\big(t^1_1,\dots, t^1_{\xi^1};\dots; t^{N-1}_1,\dots, t^{N-1}_{\xi^{N-1}}\big)$.
In the evaluation $Y({\mathfrak{gl}}_N)$-module $V(x)$, one has
\begin{gather}
{\mathbb B}_{\xi}(t) v =
\prod_{i=1}^{\xi^{N-1}}\frac{1}{t^{N-1}_i-x}
 \sum_\eta \frac{1}{\eta^{1}!} \prod_{a=1}^{N-2}
\frac{1}{(\xi^a-\eta^{a})! (\eta^{a+1}-\eta^{a})!}\nonumber\\
\hphantom{{\mathbb B}_{\xi}(t) v = }{}
\times \overline{\operatorname{Sym}}^\xi_{t}
\Bigg[X_\eta(t_{(\xi-\eta, \xi]})
Z_{\xi-\eta,\eta}(t_{[\xi-\eta]};t_{(\xi-\eta,\xi]})
\prod_{a=1}^{N-2}\prod_{i=0}^{\eta^a-1}
\frac{t^{a}_{\xi^a-i}-x+\Lambda^{a+1}}{t^{a}_{\xi^a-i}-x}
\nonumber
\\[1ex]
\hphantom{{\mathbb B}_{\xi}(t) v = }{}
\times
e_{N,N-1}^{\eta^{N-1}-\eta^{N-2}}
e_{N,N-2}^{\eta^{N-2}-\eta^{N-3}}\cdots e_{N1}^{\eta^1}
\phi\bigl({\mathbb B}^{\langle N-1\rangle}_{(\xi-\eta)^{\cdot}} (\dot t_{[\xi-\eta]})\bigr) v
\Biggr],\label{BBX}
\end{gather}
the sum being taken over all  $\eta=(\eta^{1},\dots,\eta^{N-1})\in{\mathbb Z}_{\geqslant 0}^{N-1}$ such that
$\eta^{1}\leqslant \cdots \leqslant\eta^{N-1}=\xi^{N-1}$ and $\eta^{a}\leqslant\xi^a$
for all $a=1,\dots, N-2$. Other notation is as follows: $\overline{\operatorname{Sym}}^\xi_{t}$
is defined by~\eqref{Symb}, the functions $X_\eta$ and $Z_{\xi-\eta,\eta}$ are
respectively given by formulae \eqref{X} and \eqref{Z}, $\phi$ is the first of embeddings~\eqref{embed}, and
\[
{{\mathbb B}^{\langle N-1\rangle}_{(\xi-\eta)^{\cdot}}(\dot t_{[\xi-\eta]})=
{\mathbb B}^{\langle N-1\rangle }_\zeta(s)\big|_{\zeta=(\xi-\eta)^{\cdot},\; s=\dot t_{[\xi-\eta]}}},
\]
${\mathbb B}^{\langle N-1\rangle}_\zeta(s)$ coming from \eqref{BB}.
\end{thm}

\begin{rem}
For ${N=2}$, the sum in the right side of formula \eqref{BBX} contains only one term:
${\eta=\xi}$. Moreover, ${X_\eta= Z_{\xi-\eta, \eta}= 1}$, and
${\mathbb B}^{\langle 1\rangle}_{(\xi-\eta)^{\cdot}} = 1$ by convention.
\end{rem}

\begin{cor}
\label{BCN}
Let $V$ be a ${\mathfrak{gl}}_N$-module and ${v\in V}$ a singular vector of weight $\big(\Lambda^{1},\dots,\Lambda^{N}\big)$.
Let $\xi^1,\dots,\xi^{N-1}$ be nonnegative integers and $t=\big(t^1_1,\dots, t^1_{\xi^1};\dots; t^{N-1}_1,\dots, t^{N-1}_{\xi^{N-1}}\big)$.
In the evaluation $Y({\mathfrak{gl}}_N)$-module $V(x)$, one has
\begin{gather}
{\mathbb B}_{\xi}(t)v =
\prod_{a=1}^{N-1} \prod_{i=1}^{\xi^a} \frac1{t^{a}_i-x}
\sum_m\left[ \mathop{\overleftarrow\prod}\limits_{1\leqslant b<a\leqslant N}\frac1{(m^{ab}-m^{a,b-1})!}
e_{ab}^{m^{ab}-m^{a,b-1}}\right] v\nonumber
\\[1ex]
\hphantom{{\mathbb B}_{\xi}(t)v =}{}
\times \overline{\operatorname{Sym}}^\xi_{t}\left[ \prod_{a=3}^N \prod_{b=1}^{a-2}
 \prod_{i=1}^{m^{ab}}\!\left( \frac{t^b_{i+{\widetilde m}^{ab}}-x+\Lambda^{b+1}}
{t^{b+1}_{i+{\widetilde m}^{a,b+1}}-t^b_{i+{\widetilde m}^{ab}}}
\prod_{1\leqslant j<i+{\widetilde m}^{a,b+1}}\!
\frac{t^{b+1}_j-t^b_{i+{\widetilde m}^{ab}}+1}{t^{b+1}_j-t^b_{i+{\widetilde m}^{ab}}}
\right)\right]\! .\!\!\!\!\label{BcN}
\end{gather}
Here the sum is taken over all collections of nonnegative integers $m^{ab}$,
$1\leqslant b<a\leqslant N$, such that $m^{a1}\leqslant\cdots\leqslant m^{a,a-1}$ and
$m^{a+1,a}+\cdots +m^{Na} =\xi^a$ for all $a=1,\dots,N-1$; by convention,
$m^{a0} =0$  for any $a=2,\dots, N$. Other notation is as follows:
in the ordered product the factor $e_{ab}^\circledast$ is to the left of the factor $e_{cd}^\circledast$ if $a>c$, or $a=c$ and $b>d$, $\overline{\operatorname{Sym}}^\xi_{t}$ is defined
by~\eqref{Symb}, and ${\widetilde m}^{ab}=m^{b+1,b}+\cdots +m^{a-1,b}$ for all
$1\leqslant b<a\leqslant N$, in particular, ${\widetilde m}^{a,a-1}=0$.
\end{cor}

\begin{thm}
\label{BB1}
Let $V$ be a ${\mathfrak{gl}}_N$-module and ${v\in V}$ a singular vector of weight $\big(\Lambda^{1},\dots,\Lambda^{N}\big)$.
Let $\xi^1,\dots,\xi^{N-1}$ be nonnegative integers and $t=\big(t^1_1,\dots, t^1_{\xi^1};\dots; t^{N-1}_1,\dots, t^{N-1}_{\xi^{N-1}}\big)$.
In the evaluation $Y({\mathfrak{gl}}_N)$-module $V(x)$, one has
\begin{gather}
{\mathbb B}_{\xi}(t)v = \prod_{i=1}^{\xi^1} \frac1{t^1_i-x}
 \sum_\eta  \frac1{\eta^{N-1}!}\prod_{a=2}^{N-1}
\frac1{(\xi^a-\eta^{a})!(\eta^{a-1}-\eta^{a})!}
\nonumber
\\
\hphantom{{\mathbb B}_{\xi}(t)v =}{}
\times \overline{\operatorname{Sym}}^\xi_{t} \Biggl[ Y_\eta(t_{[\eta]})
Z_{\eta,\xi-\eta}(t_{[\eta]};t_{(\eta,\xi]})\prod_{a=2}^{N-1}
\prod_{i=1}^{\eta^a}\frac{t^{a}_i-x+\Lambda^{a}}{t^{a}_i-x} \nonumber
\\
\hphantom{{\mathbb B}_{\xi}(t)v =}{}
\times e_{21}^{\eta^1-\eta^2} e_{31}^{\eta^2-\eta^3}\cdots
e_{N1}^{\eta^{N-1}}\psi\bigl({\mathbb B}^{\langle N-1\rangle }_{(\xi-\eta)^{\cdot \cdot}}
(\ddot t_{(\eta,\xi]})\bigr)v\Biggr] , \label{BBY}
\end{gather}
the sum being taken over all $\eta=(\eta^{1},\dots,\eta^{N-1})\in{\mathbb Z}_{\geqslant 0}^{N-1}$ such that
$\xi^1=\eta^{1}\geqslant\cdots\geqslant\eta^{N-1}$ and $\eta^{a}\leqslant\xi^a$
for all $a=2,\dots, N-1$. Other notation is as follows: $\overline{\operatorname{Sym}}^\xi_{t}$
is defined by~\eqref{Symb}, the functions $Y_{\eta}$ and
$Z_{\eta,\xi-\eta}$ are respectively given by formulae~\eqref{Y} and~\eqref{Z},
$\psi$ is the second of embeddings~\eqref{embed}, and
\[
{{\mathbb B}^{\langle N-1\rangle}_{(\xi-\eta)^{\cdot\cdot}}(\ddot t_{(\eta,\xi]})=
{\mathbb B}^{\langle N-1\rangle}_\zeta(s)\big|_{\zeta=(\xi-\eta)^{\cdot\cdot},\; s=\ddot t_{(\eta,\xi]}}},
\]
${\mathbb B}^{\langle N-1\rangle}_\zeta(s)$ coming from \eqref{BB}.
\end{thm}

\begin{rem}
For $N=2$ the sum in the right side of formula~\eqref{BBY} contains only one term:
$\eta=\xi$. Moreover, ${Y_\eta=Z_{\eta,\xi-\eta}=1}$, and
${\mathbb B}^{\langle 1\rangle}_{(\xi-\eta)^{\cdot\cdot}} = 1$ by convention.
\end{rem}

\begin{cor}
\label{BC1}
Let $V$ be a ${\mathfrak{gl}}_N$-module and ${v\in V}$ a singular vector of weight $\big(\Lambda^{1},\dots,\Lambda^{N}\big)$.
Let $\xi^1,\dots,\xi^{N-1}$ be nonnegative integers and $t=\big(t^1_1,\dots, t^1_{\xi^1};\dots; t^{N-1}_1,\dots, t^{N-1}_{\xi^{N-1}}\big)$.
In the evaluation $Y({\mathfrak{gl}}_N)$-module $V(x)$, one has
\begin{gather}
{\mathbb B}_{\xi}(t)v =
\prod_{a=1}^{N-1}\prod_{i=1}^{\xi^a}\frac1{t^{a}_i-x}
 \sum_m\left[ \mathop{\overrightarrow\prod}\limits_{1\leqslant b<a\leqslant N}\frac1{(m^{ab}-m^{a+1,b})!}
e_{ab}^{m^{ab}-m^{a+1,b}}\right]  v
\nonumber\\[0.5ex]
\times \overline{\operatorname{Sym}}^\xi_{t}\left[ \prod_{a=2}^{N-1}\prod_{b=1}^{a-1}
\prod_{i=0}^{m^{a+1,b}-1}\!
\left(\frac{t^{a}_{{\widehat m}^{a+1,b}-i}-x+\Lambda^{a}}
{t^{a}_{{\widehat m}^{a+1,b}-i}-t^{a-1}_{{\widehat m}^{ab}-i}}
\prod_{{\widehat m}^{ab}-i<j\leqslant\xi^{a-1}}\!
\frac{t^{a}_{{\widehat m}^{a+1,b}-i}-t^{a-1}_j+1}
{t^{a}_{{\widehat m}^{a+1,b}-i}-t^{a-1}_j}\right)\right]\!.\!\!\!\!\!\label{Bc1}
\end{gather}
Here the sum is taken over all collections of nonnegative integers $m^{ab}$,
$1\leqslant b<a\leqslant N$, such that $m^{a+1,a}\geqslant\cdots \geqslant m^{Na}$ and
$m^{a+1,1}+\cdots +m^{a+1,a}=\xi^a$ for all $a=1,\dots,N-1$; by convention,
$m^{N+1,a}=0$ for any $a=1,\dots,N$. Other notation is as follows:
in the ordered product the factor $e_{ab}^\circledast$ is to the left of the factor
$e_{cd}^\circledast$ if $b<d$, or $b=d$ and $a<c$, $\overline{\operatorname{Sym}}^\xi_{t}$ is defined by~\eqref{Symb}, and ${\widehat m}^{ab}=m^{a1}+\cdots +m^{ab}$ for all $1\leqslant b<a\leqslant N$,
in particular, ${\widehat m}^{a+1,a}=\xi^a$.
\end{cor}

\begin{thm}[\cite{TV1}]\label{BBvox}
Let $V_1$, $V_2$ be $Y({\mathfrak{gl}}_N)$-modules and ${v_1\in V_1}$, ${v_2\in V_2}$ weight
singular vectors with respect to the action of $Y({\mathfrak{gl}}_N)$. Let $\xi^1,\dots,\xi^{N-1}$ be nonnegative
integers and $t=\big(t^1_1,\dots, t^1_{\xi^1};\dots; t^{N-1}_1,\dots, t^{N-1}_{\xi^{N-1}}\big)$. Then
\begin{gather}
{\mathbb B}_{\xi}(t)(v_1\otimes v_2) =
\sum_{\eta } \prod_{a=1}^{N-1}
\frac1{(\xi^a-\eta^{a})!\eta^{a}!}
\overline{\operatorname{Sym}}^\xi_{t}\left[\prod_{a=1}^{N-2} \prod_{i=1}^{\eta^{a+1}}
\prod_{j=\eta^a +1 }^{\xi^a}
\frac{t^{a+1}_i-t^{a}_j+1}{t^{a+1}_i-t^{a}_j} \right.\nonumber
\\
\left.{}
\times
\prod_{a=1}^{N-1} \left(
\prod_{i=1}^{\eta^a} \big\langle T_{aa}(t^{a}_i) v_2\big\rangle
\prod_{j=\eta^a+1}^{\xi^a}\big\langle T_{a+1,a+1}(t^{a}_j)v_1\big\rangle
 \right) {\mathbb B}_{\eta}(t_{[\eta]}) v_1
 \otimes{\mathbb B}_{\xi-\eta}(t_{(\eta,\xi]})v_2\right],\label{BBvx}
\end{gather}
the sum being taken over all $\eta=(\eta^{1},\dots,\eta^{N-1})\in{\mathbb Z}_{\geqslant 0}^{N-1}$ such that
$\xi-\eta\in{\mathbb Z}_{\geqslant 0}^{N-1}$. In the left side we assume that ${\mathbb B}_{\xi}(t)$
acts in the $Y({\mathfrak{gl}}_N)$-module $V_1\otimes V_2$.
\end{thm}

To make the paper self-contained we will prove Theorem~\ref{BBvox} in
Section~\ref{:SBBvox}.

\begin{cor}
\label{BBvoxn}
Let $V_1,\dots, V_n$ be $Y({\mathfrak{gl}}_N)$-modules and ${v_r\in V_r}$, $r=1,\dots,n$, weight
singular vectors with respect to the action of $Y({\mathfrak{gl}}_N)$. Let $\xi^1,\dots,\xi^{N-1}$ be nonnegative
integers and $t=\big(t^1_1,\dots, t^1_{\xi^1};\dots; t^{N-1}_1,\dots, t^{N-1}_{\xi^{N-1}}\big)$. Then
\begin{gather}
 {\mathbb B}_{\xi}(t) ( v_1\otimes\cdots\otimes v_n) \nonumber\\
 \qquad{}
=  \sum_{\eta_1,\dots, \eta_{n-1}} \prod_{a=1}^{N-1}
 \prod_{r=1}^n  \frac1{(\eta^{a}_r -\eta^{a}_{r-1})!}
\overline{\operatorname{Sym}}^\xi_{t}\Biggl[\prod_{a=1}^{N-2}\prod_{r=1}^{n-1}
\prod_{i=\eta^{a+1}_{r-1}+1}^{\eta^{a+1}_r}
\prod_{j=\eta^a_r+1}^{\xi^a}
\frac{t^{a+1}_i -t^{a}_j +1}{t^{a+1}_i -t^{a}_j} \nonumber\\
\qquad\quad{}\times \prod_{a=1}^{N-1} \prod_{r=1}^n\Biggl(
\prod_{i=1}^{\eta^a_{r-1}}\big\langle T_{aa}(t^{a}_i) v_r\big\rangle
\prod_{j=\eta^a_r +1 }^{\xi^a} \big\langle T_{a+1,a+1}(t^{a}_j) v_r\big\rangle
 \Biggr) \nonumber
\\
\qquad\quad{}
\times {\mathbb B}_{\eta_1}(t_{[\eta_1]})v_1 \otimes
{\mathbb B}_{\eta_2-\eta_1}(t_{(\eta_1,\eta_2]})v_2\otimes\cdots\otimes
{\mathbb B}_{\xi-\eta_{n-1}}(t_{(\eta_{n-1},\xi]})v_n\Biggr]. \label{BBvxn}
\end{gather}
Here the sum is taken over all $\eta_1,\dots,\eta_{n-1}\in{\mathbb Z}_{\geqslant 0}^{N-1}$,
$\eta_r=\big(\eta^{1}_r,\dots,\eta^{N-1}_r\big)$, such that
$\eta_{r+1}-\eta_r\in{\mathbb Z}_{\geqslant 0}^{N-1}$ for any $r=1,\dots,n-1$, and
${\eta_0=0}$, ${\eta_n=\xi}$, by convention.
The sets $t_{[\eta_1]}$, $t_{(\eta_r,\eta_{r+1}]}$ are defined by~\eqref{teta}. In the left side we assume that ${\mathbb B}_{\xi}(t)$ acts in the $Y({\mathfrak{gl}}_N)$-module $V_1\otimes\cdots\otimes V_n$.
\end{cor}

\begin{rem}
In formulae \eqref{BBX}--\eqref{BBvxn}, the products of factorials in
the denominators of the f\/irst factors of summands are equal to the orders
of the stationary subgroups of expressions in the square brackets.
\end{rem}

\section{Proofs of Theorems \ref{BBN} and \ref{BB1}}

We prove Theorems~\ref{BBN} and \ref{BB1} by induction with respect to~$N$, assuming
that Theorem~\ref{BBvox} holds. For the base of induction, ${N=2}$, the claims
of Theorems~\ref{BBN} and~\ref{BB1} coincide with each other and reduce
to the identity
\begin{gather}
\sum_{\sigma \in S_k}  \prod_{1\leqslant i<j\leqslant k}
\frac{s_{\sigma_i}-s_{\sigma_j} -1}{s_{\sigma_i} -s_{\sigma_j}} = k!.
\label{k!}
\end{gather}
The induction step for Theorem~\ref{BBN} (resp.~\ref{BB1})
is based on Proposition~\ref{bbN} (resp.~\ref{bb1}).

Let $E^{\langle N-1\rangle}_{ab}\in\operatorname{End}\big({\mathbb C}^{N-1}\big)$ be a matrix with the only nonzero entry
equal to $1$ at the intersection of the $a$-th row and $b$-th column,
$R^{\langle N-1\rangle}(u)$ the corresponding rational $R$-matrix, cf.~\eqref{R}, and $T^{\langle N-1\rangle}_{ab}(u)$
series \eqref{Tab} for the algebra $ Y({\mathfrak{gl}}_{N-1})$. Denote by $L(x)$ a $Y({\mathfrak{gl}}_{N-1})$-module def\/ined
on the vector space ${\mathbb C}^{N-1}$ by the rule
\[
\pi(x):\ T^{\langle N-1\rangle}_{ab}(u) \mapsto \delta_{ab} + (u-x)^{-1}E^{\langle N-1\rangle}_{ba} .
\]
Denote by $\bar L(x)$ a $Y({\mathfrak{gl}}_{N-1})$-module def\/ined on the space ${\mathbb C}^{N-1}$ by the rule
\[
\varpi(x):\ T^{\langle N-1\rangle}_{ab}(u) \mapsto \delta_{ab} - (u-x)^{-1}E^{\langle N-1\rangle}_{ab} .
\]
Using $R$-matrices, the rules can be written as follows:{\samepage
\begin{gather*}
\pi(x): \ T^{\langle N-1\rangle}(u) \mapsto (u-x)^{-1}R^{\langle N-1\rangle}(u-x),
\\
\varpi(x):\ T^{\langle N-1\rangle}(u) \mapsto
(x-u)^{-1}\Bigl(\bigl(R^{\langle N-1\rangle}(x-u)\bigr)^{(21)}\Bigr)^{t_2} ,
\end{gather*}
the superscript $t_2$ standing for the matrix transposition in the second
tensor factor.}

Let ${{\bf w}}_1,\dots,{{\bf w}}_{N-1}$ be the standard basis of the space ${\mathbb C}^{N-1}$. The module $L(x)$ is
a highest weight evaluation module with highest weight $(1,0,\dots,0)$ and highest weight vector ${{\bf w}}_1$. The module $\bar L(x)$
is a highest weight evaluation module with highest weight $(0,\dots,0,-1)$ and highest weight vector ${{\bf w}}_{N-1}$.

For any ${X \in\operatorname{End}\big({\mathbb C}^{N-1}\big)}$ set  ${\nu(X)=X {{\bf w}}_1}$  and
${\bar\nu(X)=X {{\bf w}}_{N-1}}$.

Consider the maps
$\psi(x_1,\dots,x_k): Y({\mathfrak{gl}}_{N-1})\to \big({\mathbb C}^{N-1}\big)^{\otimes k} \otimes Y({\mathfrak{gl}}_N)$,
\begin{gather}
\psi(x_1,\dots,x_k)  = \big(\nu^{\otimes k}\otimes {\rm id}\big)\circ
\bigl(\pi(x_1)\otimes\cdots\otimes\pi(x_k)\otimes\psi\bigr)\circ \bigl(\Delta^{\langle N-1\rangle}\bigr)^{(k)} ,
\label{psix}
\end{gather}
and $\phi(x_1,\dots,x_k): Y({\mathfrak{gl}}_{N-1}) \to Y({\mathfrak{gl}}_N)\otimes\big({\mathbb C}^{N-1}\big)^{\otimes k} $,
\[
\phi(x_1,\dots,x_k) = \big( {\rm id}\otimes\bar\nu^{\otimes k}\big)\circ
\bigl(\phi\otimes\varpi(x_1)\otimes\cdots\otimes\varpi(x_k)\bigr)\circ \bigl(\Delta^{\langle N-1\rangle}\bigr)^{(k)},
\]
where $\psi$ and $\phi$ are embeddings \eqref{embed}, and
 ${\bigl(\Delta^{\langle N-1\rangle}\bigr)^{(k)} :
  Y({\mathfrak{gl}}_{N-1}) \to \bigl( Y({\mathfrak{gl}}_{N-1}) \bigr)^{\otimes(k+1)}}$
is the multiple coproduct.

For any element ${g\in\big({\mathbb C}^{N-1}\big)^{\otimes k}\otimes Y({\mathfrak{gl}}_N)}$ we def\/ine its components
$g^{a_1,\dots,a_k}$ by the rule
\[
g=\sum_{a_1,\dots,a_k=1}^{N-1} {{\bf w}}_{a_1}\otimes\cdots\otimes {{\bf w}}_{a_k}  \otimes g^{a_1,\dots,a_k}.
\]
A similar rule def\/ines components of elements of the tensor product
${Y({\mathfrak{gl}}_N)\otimes\big({\mathbb C}^{N-1}\big)^{\otimes k}}$.

\begin{prop}[\protect{\cite[Theorem~3.4.2]{TV1}}]
\label{bb1}
Let $\xi^1,\dots,\xi^{N-1}$ be nonnegative integers and $t=\big(t^1_1,$ $\dots, t^1_{\xi^1};\dots;t^{N-1}_1,\dots, t^{N-1}_{\xi^{N-1}}\big)$. Then
\begin{gather}
{\mathbb B}_{\xi}(t)=\sum_{a_1,\dots, a_{\xi^1}= 1}^{N-1}
T_{1,a_1+1}\big(t^1_1\big)\cdots T_{1,a_{\xi^1}+1}\big(t^1_{\xi^1}\big)
\Bigl(\psi\big(t^1_1,\dots, t^1_{\xi^1}\big)\bigl({\mathbb B}^{\langle N-1\rangle}_{\ddot\xi}(\ddot t)\bigr)\Bigr)
^{a_1,\dots,a_{\xi^1}},
\label{BB1e}
\end{gather}
cf.~\eqref{dot}.
\end{prop}

\begin{proof}
To get formula \eqref{BB1e} we use formulae \eqref{BBh} and \eqref{BB}, and compute
the trace over the f\/irst $\xi_1$ tensor factors, taking into account
the properties of the $R$-matrix~\eqref{R} described below.

Let ${{\bf v}}_1,\dots,{{\bf v}}_N$ be the standard basis of the space ${\mathbb C}^N$.
For any $a,b=1,\dots,N$, the $R$-matrix $R(u)$ preserves the subspace spanned
by the vectors ${{\bf v}}_a \otimes{{\bf v}}_b$ and ${{\bf v}}_b \otimes {{\bf v}}_a$.

Let $W$ be the image of ${{\mathbb C}^{N-1}\otimes{\mathbb C}^{N-1}}$ in ${{\mathbb C}^N \otimes{\mathbb C}^N}$ under
the embedding ${{\bf w}}_a\otimes{{\bf w}}_b\mapsto{{\bf v}}_{a+1}\otimes {{\bf v}}_{b+1}$,
$a,b=1,\dots,N-1$. The $R$-matrix $R(u)$ preserves $W$ and the restriction of $R(u)$
on $W$ coincides with the image of $R^{\langle N-1\rangle}(u)$ in $\operatorname{End}({{\mathbb C}^{N-1}\otimes{\mathbb C}^{N-1}})$
induced by the embedding.
\end{proof}

\begin{prop}
\label{bbN}
Let $\xi^1,\dots,\xi^{N-1}$ be nonnegative integers and $t=\big(t^1_1,\dots, t^1_{\xi^1};\dots; t^{N-1}_1,\dots,$ $t^{N-1}_{\xi^{N-1}}\big)$. Then
\begin{gather}
{\mathbb B}_{\xi}(t)=
\!\!\!\sum_{a_1,\dots, a_{\xi^1}=1}^{N-1}\!\!\!
T_{a_{\xi^{N-1}}+1,1}\big(t^{N-1}_{\xi^{N-1}}\big)\cdots T_{a_1+1,1}\big(t^{N-1}_1\big)
\!\Bigl(\phi\big(t^1_1,\dots, t^1_{\xi^1}\big)\bigl({\mathbb B}^{\langle N-1\rangle}_{\dot\xi}(\dot t)\bigr)
\!\Bigr)^{a_1,\dots,a_{\xi^1}}\!,\!\!\!\label{BBNe}
\end{gather}
cf.~\eqref{dot}.
\end{prop}

\begin{proof}
To get formula~\eqref{BBNe} we modify formula~\eqref{BBh} according to relation
\eqref{TTR}, use formula~\eqref{BB}, and compute the trace over the last
$\xi_{N-1}$ tensor factors, taking into account the structure of the $R$-matrix~\eqref{R}.
\end{proof}

\begin{proof}[Proof of Theorem~\ref{BB1}]
For a collection ${\boldsymbol{a}}=(a_1,\dots,a_{\xi^1})$ of positive integers let
$c^{b}({\boldsymbol{a}})=\#
\{r\vert a_r\geqslant b\}$, and
$c({\boldsymbol{a}})=\bigl(c^{1}({\boldsymbol{a}}),\dots, c^{N-1}({\boldsymbol{a}})\bigr)$.

To obtain formula \eqref{BBY} we apply both sides of formula \eqref{BB1e} to the
singular vector $v$ in the evaluation module $V(x)$ over $Y({\mathfrak{gl}}_N)$. In this case, $T_{1a}(u)$
acts as $(u-x)^{-1}e_{a1}$ and we have
\begin{gather}
{\mathbb B}_{\xi}(t)v = \prod_{i=1}^{\xi^1}\frac1{t^1_i -x}\nonumber\\
\hphantom{{\mathbb B}_{\xi}(t)v =}{}\times
\sum_\eta e_{21}^{\eta^1-\eta^2}
e_{31}^{\eta^2-\eta^3}\cdots e_{N1}^{\eta^{N-1}}
 \sum_{\substack{a_1,\dots, a_{\xi^1}=1 \\ c({\boldsymbol{a}})=\eta}}^{N-1}
 \Bigl(\psi\big(t^1_1,\dots, t^1_{\xi^1}\big)\bigl({\mathbb B}^{\langle N-1\rangle}_{\ddot\xi}(\ddot t)\bigr)
\Bigr)^{a_1,\dots, a_{\xi^1}} v,\label{BBci}
\end{gather}
the f\/irst sum being taken over all $\eta=\big(\eta^{1},\dots,\eta^{N-1}\big)\in{\mathbb Z}_{\geqslant 0}^{N-1}$ such that
$\xi^1=\eta^{1}\geqslant\cdots\geqslant\eta^{N-1}$.

Let $^\psi V(x)$ be the $Y({\mathfrak{gl}}_{N-1})$-module obtained by pulling $V(x)$ back
through the embedding $\psi$. Then
$\psi\big(t^1_1,\dots, t^1_{\xi^1}\big)\bigl({\mathbb B}^{\langle N-1\rangle}_{\ddot\xi}(\ddot t)\bigr) v$ is the weight function
associated with the vector ${{{\bf w}}_1\otimes\cdots\otimes {{\bf w}}_1 \otimes v}$ in the $Y({\mathfrak{gl}}_{N-1})$-module
$L\big(t^1_1\big)\otimes\cdots\otimes L\big(t^1_{\xi^1}\big)\otimes {}^\psi V(x)$. We use Theorem~\ref{BBvox} to write
$\psi\big(t^1_1,\dots, t^1_{\xi^1}\big)\bigl({\mathbb B}^{\langle N-1\rangle}_{\ddot\xi}(\ddot t)\bigr) v$ as a sum of tensor
products of weight functions in the tensor factors, that is, as a sum of the
following expressions:
\[
\pi\big(t^1_1\big)\bigl({\mathbb B}^{\langle N-1\rangle}_{\zeta_1}(s_1)\bigr) {{\bf w}}_1\otimes\cdots\otimes
\pi\big(t^1_{\xi^1}\big)\bigl({\mathbb B}^{\langle N-1\rangle}_{\zeta_{\xi^1}}(s_{\xi^1})\bigr) {{\bf w}}_1
\otimes \psi\bigl({\mathbb B}^{\langle N-1\rangle}_{\zeta_0}(s_0)\bigr) v ,
\]
where $\zeta_0,\dots, \zeta_{\xi^1}$, $s_0,\dots, s_{\xi^1}$ are suitable parameters,
and employ Corollary~\ref{BC1}, valid by the induction assumption, to calculate
the weight functions $\pi(t^1_j)\bigl({\mathbb B}^{\langle N-1\rangle}_{\zeta_j}(s_j)\bigr) {{\bf w}}_1$ in
the modules~$L(t^1_j)$. As a result, we get formula~\eqref{psihBB},
see~Lemma~\ref{psihbb} below.

Observe that in the module $L(x)$ one has
$\langle T_{11}(u) {{\bf w}}_1\rangle = 1+(u-x)^{-1}$ and
$\langle T_{aa}(u) {{\bf w}}_1\rangle= 1$ for all $a=2,\dots, N$.
The weight function $\pi(x)\bigl({\mathbb B}^{\langle N-1\rangle}_\zeta(s)\bigr) {{\bf w}}_1$ equals zero
unless $\zeta=(1,\dots,1,0,\dots,0)$ (it can be no units or zeros in the sequence).
If $\zeta^1=\cdots=\zeta^r =1$ and  $\zeta^{r+1}=\cdots =\zeta^{N-1}=0$,
then $s=(s^1_1,\dots, s^r_1)$ and
\[
\pi(x)\bigl({\mathbb B}^{\langle N-1\rangle}_\zeta(s)\bigr) {{\bf w}}_1 = \frac{e_{r+1,1} {{\bf w}}_1}
{(s^1_1 -x) (s^2_1 -s^1_1)\cdots(s^r_1 -s^{r-1}_1)} .
\]

Fix $\eta=\big(\eta^{1},\dots,\eta^{N-1}\big)\in{\mathbb Z}_{\geqslant 0}^{N-1}$ such that $\eta^{1}\geqslant\cdots \geqslant\eta^{N-1}$.
Consider a collection ${\boldsymbol{l}}$ of integers $l^{a}_i$,
$a=1,\dots,N-2$, $i=1,\dots,\eta^{a+1}$, such that
$1\leqslant l^{a}_1<\cdots <l^{a}_{\eta^{a+1}}\leqslant\eta^{a}$ for all $a=1,\dots,N-2$.
Introduce a function $F_{{\boldsymbol{l}}}(s)$ of the variables $s^1_1,\dots,s^1_{\eta^1};\dots; s^{N-1}_1,\dots,s^{N-1}_{\eta^{N-1}}$:
\begin{gather}
F_{{\boldsymbol{l}}}(s)=\prod_{a=1}^{N-2} \prod_{i=1}^{\eta^{a+1}}\left(
\frac1{s^{a+1}_i -s^a_{l^a_i}} \prod_{l^a_i<j\leqslant\eta^{a}}
\frac{s^{a+1}_i-s^a_j+1}{s^{a+1}_i-s^a_j}\right).
\label{F}
\end{gather}
There is a bijection between collections ${\boldsymbol{l}}$ and sequences of integers
${\boldsymbol{a}}=(a_1,\dots, a_{\eta^1})$ such that $1\leqslant a_i\leqslant N-1$ for all $i=1,\dots,\eta^{1}$,
and $c({\boldsymbol{a}})=\eta$. It is established as follows. Def\/ine numbers~$p^{a}_i$
by the rule: $p^{1}_i=l^{1}_i$, $i=1,\dots,\eta^{2}$, and
$p^{a}_i=p^{a-1}_{l^a_i}$, $a=2,\dots,N-2$, $i=1,\dots,\eta^{a+1}$.
Then the sequence ${\boldsymbol{a}}$ is uniquely determined by the requirement that $a_i>b$
if\/f $i\in\big\{p^{b}_1,\dots, p^{b}_{\eta^{b+1}}\big\}$, for all
$i=1,\dots,\eta^{1}$. We will write ${\boldsymbol{a}}({\boldsymbol{l}})$ for the result
of this mapping.

Summarizing, we get the following statement.

\begin{lem}
\label{psihbb}
Let $\eta=(\eta^{1},\dots,\eta^{N-1})\in{\mathbb Z}_{\geqslant 0}^{N-1}$ be such that
$\xi^1=\eta^{1}\geqslant\cdots \geqslant\eta^{N-1}$. Let ${\boldsymbol{l}}$ be a collection
of integers as described above, and ${\boldsymbol{a}}({\boldsymbol{l}})=(a_1,\dots, a_{\xi^1})$. Then
\begin{gather}
\Bigl(\psi\big( t^1_1,\dots, t^1_{\xi^1}\big)\bigl({\mathbb B}^{\langle N-1\rangle}_{\ddot\xi}(\ddot t)\bigr)
\Bigr)^{a_1,\dots, a_{\xi^1}}v\label{psihBB}\\[0.5ex]
=\prod_{b=2}^{N-1}\frac1{(\xi^b-\eta^{b})!}\,
\overline{\operatorname{Sym}}^{\ddot\xi}_{\ddot t} \left[F_{{\boldsymbol{l}}}(t_{[\eta]})
Z_{\eta,\xi-\eta}(t_{[\eta]};t_{(\eta,\xi]})
\prod_{b=2}^{N-1}\prod_{i=1}^{\eta^b}
\frac{t^b_i-x+\Lambda^{b}}{t^b_i-x}
 \psi\bigl({\mathbb B}^{\langle N-1\rangle}_{(\xi-\eta)^{\cdot\cdot}}
(\ddot t_{(\eta,\xi]})\bigr) v \right],\nonumber
\end{gather}
cf.~\eqref{Z} for $Z_{\eta,\xi-\eta}(t_{[\eta]};t_{(\eta,\xi]})$.
\end{lem}

Comparing the expressions under $\overline{\operatorname{Sym}}$ in formulae \eqref{psihBB} and
\eqref{BBY}, and taking into account that the product
\[
Z_{\eta,\xi-\eta}(t_{[\eta]};t_{(\eta,\xi]})
\prod_{b=2}^{N-1}\prod_{i=1}^{\eta^b}\big(t^b_i-x+\Lambda^{b}\big)
\psi\bigl({\mathbb B}^{\langle N-1\rangle}_{(\xi-\eta)^{\cdot\cdot}}
(\ddot t_{(\eta,\xi]})\bigr) v
\]
is invariant with respect to the action of the groups ${S_{\eta^1}\times\cdots\times S_{\eta^{N-1}}}$ and $S_{\xi^1-\eta^1}\times\cdots\times  S_{\xi^{N-1}-\eta^{N-1}}$ permuting
respectively the variables $t^1_1,\dots, t^1_{\eta^1};\dots; t^{N-1}_1,\dots, t^{N-1}_{\eta^{N-1}}$ and $t^1_{\eta_1+1},\dots$, $t^1_{\xi^1};\dots;t^{N-1}_{\eta_{N-1}+1},\dots,t^{N-1}_{\xi^{N-1}}$, one can see that formula~\eqref{BBY}
follows from formula~\eqref{BBci} and Lem\-ma~\ref{FY} below. Theorem~\ref{BB1} is proved.
\end{proof}

\begin{lem}
\label{FY}
Let $\eta=(\eta^{1},\dots,\eta^{N-1})\in{\mathbb Z}_{\geqslant 0}^{N-1}$ such that $\eta^{1}\geqslant \cdots \geqslant\eta^{N-1}$.
Let $s=s^1_1,\dots,s^1_{\eta^1};\dots;$ $s^{N-1}_1,\dots,s^{N-1}_{\eta^{N-1}}$. Then
\begin{gather}
\frac1{\eta^{N-1}!}\prod_{a=1}^{N-2}\frac1{(\eta^{a}-\eta^{a+1})!}
\overline{\operatorname{Sym}}^\eta_s\bigl(Y_\eta(s)\bigr)=\sum_{{\boldsymbol{l}}}
\overline{\operatorname{Sym}}^{\ddot\eta}_{\ddot s}\bigl(F_{{\boldsymbol{l}}}(s)\bigr),
\label{FlY}
\end{gather}
cf.~\eqref{Y} for $Y_\eta(s)$. The sum is taken over all collections ${\boldsymbol{l}}$
of integers $l^{a}_i$, $a=1,\dots,N-2$, $i=1,\dots,\eta^{a+1}$, such that
$1\leqslant l^{a}_1<\cdots <l^{a}_{\eta^{a+1}}\leqslant \eta^{a}$ for all $a=1,\dots,N-2$.
\end{lem}

\begin{proof}
It is convenient to rewrite formula \eqref{Y} in the form similar to~\eqref{F}:
\[
Y_\eta(s)=\prod_{a=1}^{N-2}\prod_{i=1}^{\eta^{a+1}}
\left(\frac1{s^{a+1}_i-s^a_{i+\eta^a-\eta^{a+1}}}
\prod_{i<j\leqslant\eta^{a+1}}
\frac{s^{a+1}_i-s^a_{j+\eta^a-\eta^{a+1}}+1}
{s^{a+1}_i-s^a_{j+\eta^a-\eta^{a+1}}}\right).
\]
For positive integers $p$, $r$ such that $p\leqslant r$, consider a function
\begin{gather*}
G_{p,r}(y_1,\dots,y_p;z_1,\dots,z_r)\\
\qquad {} =
\frac1{(r-p)!}\overline{\operatorname{Sym}}^r_{z_1,\dots,z_r}\left[
\prod_{i=1}^p \left( \frac1{y_i -z_{i+r-p}}\prod_{i<j\leqslant p}
\frac{y_i-z_{j+r-p}+1}{y_i-z_{j+r-p}}\right)\right].
\end{gather*}
It is a manifestly symmetric function of $z_1,\dots,z_r$, and it is a symmetric function of $y_1,\dots,y_p$
by Lemma~\ref{Gyz}.

To prove Lemma \ref{FY} we will show that the expressions in both sides of
formula \eqref{FlY} are equal to
\[
\prod_{a=1}^{N-2} G_{\eta^{a+1},\eta^a}
\big(s^{a+1}_1,\dots, s^{a+1}_{\eta^{a+1}};s^a_1,\dots,s^a_{\eta^a}\big).
\]
The proof is by induction with respect to $N$. The base of induction is $N=3$.
In this case the claim follows from Lemma~\ref{Gyz} and identity~\eqref{k!}.
The induction step for the left side of \eqref{FlY} is as follows:
\begin{gather*}
\overline{\operatorname{Sym}}^\eta_s \bigl(Y_\eta(s)\bigr)=
\overline{\operatorname{Sym}}^{\ddot\eta}_{\ddot s}\Bigl(
\overline{\operatorname{Sym}}^{\eta^1}_{s^1_1,\dots,s^1_{\eta^1}} \bigl(Y_\eta(s)\bigr) \Bigr)
\\
\hphantom{\overline{\operatorname{Sym}}^\eta_s \bigl(Y_\eta(s)\bigr)}{}
=\overline{\operatorname{Sym}}^{\ddot\eta}_{\ddot s} \Biggl[ G_{\eta^2,\eta^1}
\big(s^2_1,\dots,s^2_{\eta^2};s^1_1,\dots,s^1_{\eta^1}\big)
\\
\left.\hphantom{\overline{\operatorname{Sym}}^\eta_s \bigl(Y_\eta(s)\bigr)=}{}
\times \prod_{a=2}^{N-2}\prod_{i=1}^{\eta^{a+1}}
\left(\frac1{s^{a+1}_i-s^a_{i+\eta^a-\eta^{a+1}}}
\prod_{i<j\leqslant\eta^{a+1}}
\frac{s^{a+1}_i-s^a_{j+\eta^a-\eta^{a+1}}+1}
{s^{a+1}_i-s^a_{j+\eta^a-\eta^{a+1}}}\right)\right]
\\
\hphantom{\overline{\operatorname{Sym}}^\eta_s \bigl(Y_\eta(s)\bigr)}{}
=G_{\eta^2,\eta^1}
\big(s^2_1,\dots, s^2_{\eta^2};s^1_1,\dots, s^1_{\eta^1})
\\
\hphantom{\overline{\operatorname{Sym}}^\eta_s \bigl(Y_\eta(s)\bigr)=}{}
\times
\overline{\operatorname{Sym}}^{\ddot\eta}_{\ddot s}\left[
\prod_{a=2}^{N-2}\prod_{i=1}^{\eta^{a+1}}
\left(\frac1{s^{a+1}_i-s^a_{i+\eta^a-\eta^{a+1}}}
\prod_{i<j\leqslant\eta^{a+1}}
\frac{s^{a+1}_i-s^a_{j+\eta^a-\eta^{a+1}}+1}
{s^{a+1}_i-s^a_{j+\eta^a-\eta^{a+1}}}\right)\right]
\\
\hphantom{\overline{\operatorname{Sym}}^\eta_s \bigl(Y_\eta(s)\bigr)}{}
=\prod_{a=1}^{N-2}G_{\eta^{a+1},\eta^a}
\big(s^{a+1}_1,\dots, s^{a+1}_{\eta^{a+1}};s^a_1,\dots,s^a_{\eta^a}\big).
\end{gather*}
In the last two equalities we use the fact that
$G_{\eta^2,\eta^1}\big(s^2_1,\dots,s^2_{\eta^2};s^1_1,\dots,s^1_{\eta^1}\big)$
is symmetric with respect to $s^2_1,\dots,s^2_{\eta^2}$, and the induction assumption.

The idea of the induction step for the right side of~\eqref{FlY} is similar.
First, one should symmetrize $F_{{\boldsymbol{l}}}(s)$ with respect to the variables
$s^{N-1}_1,\dots,s^{N-1}_{\eta^{N-1}}$ and sum up over all possible collections
$l^{N-2}_1,\dots,l^{N-2}_{\eta^{N-1}}$, and then use Lemma~\ref{Gyz}.
We leave details to a reader.
\end{proof}

\begin{lem}
\label{Gyz}
\[
G_{p,r}(y_1,\dots,y_p;z_1,\dots,z_r)=\sum_{{\boldsymbol{d}}} \overline{\operatorname{Sym}}^p_{y_1,\dots,y_p}
\left[\prod_{i=1}^p\left(\frac1{y_i-z_{d_i}}
\prod_{d_i<j\leqslant r} \frac{y_i-z_j+1}{y_i-z_j}\right)\right],
\]
the sum being taken over all $p$-tuples ${\boldsymbol{d}}=(d_1,\dots,d_p)$ such that
${1\leqslant d_1<\dots <d_p\leqslant r}$.
\end{lem}

The proof is given at the end of Section \ref{:Twtf}.

\begin{lem}
\label{Gyz1}
Let $p$, $r$ be positive integers such that $p\leqslant r$. Then
\begin{gather*}
\frac1{(r-p)!} \overline{\operatorname{Sym}}^r_{z_1,\dots,z_r}\left[
\prod_{i=1}^p\left( \frac1{y_i-z_i}
\prod_{1\leqslant j<i} \frac{y_i-z_j+1}{y_i-z_j}\right)\right]
\\
\qquad{} =
\sum_{{\boldsymbol{d}}} \overline{\operatorname{Sym}}^p_{y_1,\dots,y_p}
\left[\prod_{i=1}^p \left(\frac1{y_i-z_{d_i}}
\prod_{1\leqslant j<d_i}\frac{y_i-z_j+1}{y_i-z_j}\right)\right],
\end{gather*}
the sum being taken over all $p$-tuples ${\boldsymbol{d}}=(d_1,\dots,d_p)$ such that
$1\leqslant d_1<\cdots <d_p\leqslant r$.
\end{lem}

\begin{proof}
The statement follows from Lemma~\ref{Gyz} by the change of variables
$y_i\to-y_{p-i}$, $z_j\to-z_{r-j}$, and a suitable change
of summation indices.
\end{proof}
\begin{proof}[Proof of Theorem~\ref{BBN}]
The proof is similar to the proof of Theorem~\ref{BB1}, mutatis mutandis.
In particular, Lemma~\ref{Gyz} is replaced by Lemma~\ref{Gyz1}.
\end{proof}

\section{Proof of Theorem~\ref{BBvox}}
\label{:SBBvox}
The theorem is proved by induction with respect to $N$. The base of induction,
the ${N=2}$ case, follows from Proposition~\ref{DlB}.
The induction step is provided by Proposition~\ref{Step}.

Let $P^{\langle N-1\rangle}=\sum\limits_{a,b=1}^{N-1} E^{\langle N-1\rangle}_{ab} \otimes E^{\langle N-1\rangle}_{ba}$ be
the f\/lip matrix, and $R^{\langle N-1\rangle}(u)=u+P^{\langle N-1\rangle}$ be the $R$-matrix for the Yangian
$ Y({\mathfrak{gl}}_{N-1})$.

In this section we regard $T(u)$ as an $N\times N$ matrix over the algebra
$Y({\mathfrak{gl}}_N)[u^{-1}]$ with entries $T_{ab}(u)$, $a,b=1,\dots,N$. Let
\begin{gather}
A(u)=T_{11}(u),\qquad B(u)=\bigl(T_{12}(u),\dots, T_{1N}(u)\bigr),
\qquad D(u)=\bigl(T_{ij}(u)\bigr)_{i,j=2}^N,
\label{ABD}
\end{gather}
be the submatrices of $T(u)$.
Set ${{\overline{R}}(u)=u^{-1}R^{\langle N-1\rangle}(u)}$. Formulae~\eqref{RTT} and~\eqref{R} imply
the following commutation relations for
$A(u)$, $B(u)$ and $D(u)$:
\begin{gather}
A(u) A(t)=A(t)A(u),
\label{AA}
\\
B^{[1]}(u) B^{[2]}(t)=\frac{u-t}{u-t+1}B^{[2]}(t)B^{[1]}(u) {\overline{R}}^{(12)}(u-t),
\label{BBR}
\\
A(u)B(t)=\frac{u-t-1}{u-t}B(t)A(u)+\frac1{u-t}B(u)A(t),
\label{AB}
\\
D^{(1)}(u)B^{[2]}(t) =\frac{u-t+1}{u-t}B^{[2]} (t) D^{(1)} (u){\overline{R}}^{(12)}(u-t)-\frac1{u-t}B^{[1]}(u)D^{(2)}(t),
\label{DB}
\\
{\overline{R}}^{(12)}(u-t) D^{(1)}(u) D^{(2)}(t)=D^{(2)}(t)D^{(1)}(u){\overline{R}}^{(12)}(u-t).
\label{DD}
\end{gather}
In this section we use superscripts to deal with tensor products of matrices,
writing parentheses for square matrices and brackets for the row matrix~$B$.

Set ${\check R}(u)=(u+1)^{-1}P^{\langle N-1\rangle}R^{\langle N-1\rangle}(u)$. For an expression $f(u_1,\dots,u_k)$
with matrix coef\/f\/i\-cients and a simple transposition $(i,i+1)$, $i=1,\dots,k-1$, set
\begin{gather}
{}^{(i,i+1)}\!f(u_1,\dots,u_k)=f(u_1,\dots, u_{i-1},u_{i+1},u_i,u_{i+2},\dots, u_k)
{\check R}^{(i,i+1)}(u_i-u_{i+1}),
\label{Ract}
\end{gather}
if the product in the right side makes sense.
The matrix ${\check R}(u)$ has the properties
${\check R}(u){\check R}(-u)=1$ and
\[
{\check R}^{(12)}(u-v){\check R}^{(23)}(u){\check R}^{(12)}(v)=
{\check R}^{(23)}(v) {\check R}^{(12)}(u){\check R}^{(23)}(u-v),
\]
cf.~\eqref{YB}. This yields the following lemma.

\begin{lem}
\label{RSym}
Formula \eqref{Ract} extends to the action of the symmetric group $S_k$ on expressions
$f(u_1,\dots,u_k)$ with appropriate matrix coefficients:
${f\mapsto {}^{\sigma}\!f}$, ${\sigma\in S_k}$.
\end{lem}

By formula \eqref{BBR} the expression $B^{[1]}(u_1)\cdots B^{[k]}(u_k)$ is invariant
under the action \eqref{Ract} of the symmetric group $S_k$.

For an expression $f(u_1,\dots,u_k)$ with suitable matrix coef\/f\/icients, set
\[
{}^{R}\!\operatorname{Sym}_{u_1,\dots,u_k}^{(1,\dots,k)} f(u_1,\dots,u_k)=\sum_{\sigma\in S_k} {}^{\sigma}\!f(u_1,\dots,u_k).
\]

\begin{prop}\label{ADB}
\begin{gather}
 A(u) B^{[1]}(u_1)\cdots B^{[k]}(u_k)=
\prod_{i=1}^k \frac{u-u_i -1}{u-u_i}  B^{[1]}(u_1)\cdots B^{[k]}(u_k) A(u)
\label{ABB}
\\
\qquad{}+\frac1{(k-1)!}\,
{}^{R}\!\operatorname{Sym}_{u_1,\dots,u_k}^{(1,\dots,k)}\left(\frac1{u-u_1}\prod_{i=2}^k
\frac{u_1-u_i-1}{u_1-u_i} B^{[1]}(u)B^{[2]}(u_2)\cdots B^{[k]}(u_k)A(u_1)
\right) ,\nonumber
\\[1ex]
  D^{(0)}(u) B^{[1]}(u_1)\cdots B^{[k]}(u_k)
\nonumber
\\
\qquad{} =  \prod_{i=1}^k \frac{u-u_i+1}{u-u_i}
\ B^{[1]}(u_1)\cdots B^{[k]}(u_k) D^{(0)}(u)
{\overline{R}}^{(0k)}(u-u_k)\cdots {\overline{R}}^{(01)}(u-u_1)\nonumber
\\
\qquad\quad{}-\frac1{(k-1)!} \, {}^{R}\!\operatorname{Sym}_{t_1,\dots,t_k}^{(1,\dots,k)}\Biggl(\frac1{u-u_1}
\prod_{i=2}^k\frac{u_1-u_i+1}{u_1-u_i} \nonumber
\\
\qquad\quad{} \times B^{[0]}(u) B^{[2]}(u_2)\cdots B^{[k]}(u_k)D^{(1)}(u_1)
{\overline{R}}^{(1k)}(u_1-u_k)\cdots {\overline{R}}^{(12)}(u_1-u_2) \Biggr).\label{DBB}
\end{gather}
In the second formula the tensor factors are counted by $0,\dots,k$.
\end{prop}

\begin{proof}
The statement follows from relations \eqref{BBR}--\eqref{DB} by induction with respect to~$k$.
We apply formula \eqref{AB} or \eqref{DB} to the product of the f\/irst factors in the left side
and then use the induction assumption.
\end{proof}

\begin{rem}
Formulae \eqref{ABB} and \eqref{DBB} have the following structure. The f\/irst term
in the right side comes from repeated usage of the f\/irst term in the right side of the respective
relation~\eqref{AB} or~\eqref{DB}. The second term, involving symmetrization,
is ef\/fectively determined by the fact that the whole expression in the right side is
regular at $u=u_i$ for any $i={1,\dots,k}$, and is invariant with respect to action~\eqref{Ract}
of the symmetric group~$S_k$. The symmetrized expression is obtained by applying once
the second term in the right side of the relevant relation~\eqref{AB} or~\eqref{DB} followed
by repeated usage of the f\/irst term of the respective relation.
\end{rem}

Let $\Delta$ be coproduct \eqref{Dl} for the Yangian $Y({\mathfrak{gl}}_N)$. For a matrix
${F=(F_{ij})}$ over $Y({\mathfrak{gl}}_N)$, denote by ${\Delta(F)=\bigl(\Delta(F_{ij})\bigr)}$
the corresponding matrix over ${Y({\mathfrak{gl}}_N)\otimes Y({\mathfrak{gl}}_N)}$.

We will use subscripts in braces to describe the embeddings $Y({\mathfrak{gl}}_N)\to Y({\mathfrak{gl}}_N)\otimes Y({\mathfrak{gl}}_N)$
as one of the tensor factors: ${X_{\{1\}}=X\otimes 1}$, ${X_{\{2\}} =1\otimes X}$,
${X \in Y({\mathfrak{gl}}_N)}$. For a matrix $F$ over $Y({\mathfrak{gl}}_N)$, we apply the embeddings entrywise,
writing $F_{\{1\}}$, $F_{\{2\}}$ for the corresponding matrices over ${Y({\mathfrak{gl}}_N)\otimes Y({\mathfrak{gl}}_N)}$.

\begin{prop}\label{DlB}
We have
\begin{gather}
\Delta\bigl(B^{[1]}(t_1)\cdots B^{[k]}(t_k)\bigr)=
\sum_{l=0}^k \frac1{l! (k-l)!}\,{}^{R}\!\operatorname{Sym}_{t_1,\dots,t_k}^{(1,\dots,k)}\Biggl(
\prod_{1\leqslant i<j\leqslant k} \frac{t_i-t_j-1}{t_i-t_j}
\nonumber\\[0.5ex]
\hphantom{\Delta\bigl(B^{[1]}(t_1)\cdots B^{[k]}(t_k)\bigr)=}{}
\times  B^{[1]}_{\{1\}}(t_1)\cdots B^{(l)}_{\{1\}}(t_l) B^{(l+1)}_{\{2\}}(t_{l+1})\cdots
B^{[k]}_{\{2\}}(t_k)\nonumber
\\
\hphantom{\Delta\bigl(B^{[1]}(t_1)\cdots B^{[k]}(t_k)\bigr)=}{}
\times
D^{(l+1)}_{\{1\}}(t_{l+1})\cdots D^{(k)}_{\{1\}}(t_k) A_{\{2\}}(t_1)\cdots A_{\{2\}}(t_l)\Biggr).\label{DlBk}
\end{gather}
\end{prop}

\begin{proof}
The statement is proved by induction with respect to $k$. Writing the left side as
\[
\Delta\bigl(B^{[1]}(u_1)\bigr) \Delta\bigl(B^{[2]}(u_2)\cdots B^{[k]}(u_k)\bigr),
\]
we expand the f\/irst factor according to \eqref{Dl}:
\[
\Delta\bigl(B^{[1]}(u_1)\bigr)=
B^{[1]}_{\{1\}}(u_1) A_{\{2\}}(u_1)+B^{[1]}_{\{2\}}(u_1)D^{(1)}_{\{1\}}(u_1),
\]
and apply the induction assumption to expand the second one. Then we use
Proposition~\ref{ADB} to transform the obtained expression to the right side of~\eqref{DlBk}.
\end{proof}

Regard vectors in the space ${\mathbb C}^{N-1}$ as $(N-1)\times 1$ matrices.
Formula \eqref{BB1e} from Proposition~\ref{bb1} can be written as follows:
\begin{gather}
{\mathbb B}_{\xi}(t)=B^{[1]}\big(t^1_1\big)\cdots B^{[\xi^1]}\big(t^1_{\xi^1}\big)
 \psi\big(t^1_1,\dots, t^1_{\xi^1}\big)\bigl( {\mathbb B}^{\langle N-1\rangle}_{\ddot\xi}(\ddot t )\bigr) .
\label{BB1i}
\end{gather}
For nonnegative integers $k$, $l$ such that $k\geqslant l$, def\/ine an embedding
\begin{gather}
\widehat\psi_{l}(u_1,\dots,u_k) : \  Y({\mathfrak{gl}}_{N-1}) \to \big({\mathbb C}^{N-1}\big)^{\otimes k} \otimes Y({\mathfrak{gl}}_N)\otimes Y({\mathfrak{gl}}_N) ,
\label{psih}
\\
\widehat\psi_{l}(u_1,\dots,u_k) =
\vartheta_l\circ \big(\nu^{ \otimes l} \otimes {\rm id}\otimes\nu^{ \otimes(k-l)} \otimes {\rm id} \big) \nonumber
\\
\hphantom{\widehat\psi_{l}(u_1,\dots,u_k) = }{}
\circ \bigl(\pi(u_1)\otimes\cdots\otimes\pi(u_l)\otimes\psi\otimes
\pi(u_{l+1})\otimes\cdots\otimes\pi(u_k)\otimes\psi\bigr)\circ \bigl(\Delta^{\langle N-1\rangle}\bigr)^{(k+1)},\nonumber
\end{gather}
where
\[
\vartheta_l: \ \big({\mathbb C}^{N-1}\big)^{\otimes l} \otimes Y({\mathfrak{gl}}_N)\otimes\big({\mathbb C}^{N-1}\big)^{\otimes(k-l)} \otimes  Y({\mathfrak{gl}}_N) \to
\big({\mathbb C}^{N-1}\big)^{\otimes k} \otimes  Y({\mathfrak{gl}}_N)\otimes  Y({\mathfrak{gl}}_N)
\]
is given by the rule
\[
\vartheta_l({\boldsymbol{x}}\otimes X_1 \otimes {\boldsymbol{y}} \otimes X_2) = {\boldsymbol{x}}\otimes{\boldsymbol{y}} \otimes X_1 \otimes X_2,
\]
for ${\boldsymbol{x}}\in\big({\mathbb C}^{N-1}\big)^{\otimes l}$, ${\boldsymbol{y}}\in\big({\mathbb C}^{N-1}\big)^{\otimes(k-l)}$,
$X_1,X_2\in Y({\mathfrak{gl}}_N)$, and
\[
\bigl(\Delta^{\langle N-1\rangle}\bigr)^{(k+1)} : \  Y({\mathfrak{gl}}_{N-1}) \to \bigl( Y({\mathfrak{gl}}_{N-1})\bigr)^{\otimes(k+2)}
\]
is the multiple coproduct.

Let $v_1$, $v_2$ be weight singular vectors with respect to the action of $ Y({\mathfrak{gl}}_N)$.

\begin{lem}
\label{DlN}
For any $X\in Y({\mathfrak{gl}}_{N-1})$ we have
\[
\Delta\bigl(\widehat\psi(u_1,\dots,u_k)(X)\bigr)(v_1\otimes v_2)
=\widehat\psi_{k}(u_1,\dots,u_k)(X)\bigr)(v_1\otimes v_2).
\]
\end{lem}

\begin{proof}
Recall that $\Delta^{\langle N-1\rangle}$ denotes the coproduct for the Yangian $ Y({\mathfrak{gl}}_{N-1})$.
Let $ Y_{\times}({\mathfrak{gl}}_N)$ be the left ideal in $ Y({\mathfrak{gl}}_N)$ generated by the coef\/f\/icients of the series
$T_{21}(u),\dots, T_{N1}(u)$. It follows from relations~\eqref{Dl} and~\eqref{Tabcd} that
\[
\Delta\bigl(\psi(X)\bigr)-(\psi\otimes\psi)\bigl(\Delta^{\langle N-1\rangle}(X)\bigr)
\in Y({\mathfrak{gl}}_N)\otimes Y_{\times}({\mathfrak{gl}}_N)
\]
for any $X\in Y({\mathfrak{gl}}_{N-1})$. Therefore,
\begin{gather}
\Delta\bigl(\psi(X)\bigr)(v_1\otimes v_2)=
(\psi\otimes\psi) \bigl(\Delta^{\langle N-1\rangle}(X)\bigr) (v_1 \otimes v_2)
\label{Dlpsi}
\end{gather}
because $v_2$ is a weight singular vector. The lemma follows from
formulae~\eqref{psix}, \eqref{psih} \linebreak and~\eqref{Dlpsi}.
\end{proof}

\begin{lem}\label{ADX}
For any $X \in Y({\mathfrak{gl}}_{N-1})$ we have
\begin{gather*}
\Bigl(D^{(l+1)}_{\{1\}} (u_{l+1})\cdots D^{(k)}_{\{1\}}(u_k)
A_{\{2\}}(t_1)\cdots A_{\{2\}}(t_l) \widehat\psi_{k}(t_1,\dots,t_k)(X) \Bigr) (v_1 \otimes v_2)
\\
\qquad {} =\prod_{i=1}^l \big\langle T_{22}\big(t^1_i\big)v_1\big\rangle
\prod_{j=l+1}^k\big\langle T_{11}\big(t^1_j\big) v_2\big\rangle \Bigl(
  \widehat\psi_{l}(t_1,\dots,t_k)(X)I^{(1)} \cdots I^{(k)}\Bigr)(v_1 \otimes v_2).
\end{gather*}
\end{lem}

\begin{proof}
Recall that $D(u)=({\rm id}\otimes\psi)\bigl(T^{\langle N-1\rangle}(u)\bigr)$ and
${\overline{R}}(u-u_i)=\bigl({\rm id}\otimes\pi(u_i)\bigr)\bigl(T^{\langle N-1\rangle}(u)\bigr)$.
Then according to relation~\eqref{DD}, for any ${X\in Y({\mathfrak{gl}}_{N-1})}$ we have
\[
D(u_i) \bigl(\psi\otimes\pi(u_i)\bigr)\bigl(\Delta^{\langle N-1\rangle}(X)\bigr)=
\bigl(\pi(u_i)\otimes\psi\bigr)\bigl(\Delta^{\langle N-1\rangle}(X)\bigr) D(u_i).
\]
In addition, remind that ${D(u)({{\bf w}}_1\otimes v_1)=
{{\bf w}}_1\otimes T_{22}(u)v_1=\big\langle  T_{22}(u) v_1\big\rangle ({{\bf w}}_1\otimes v_1)}$,
because $v_1$ is a weight singular vector. Therefore,
\begin{gather*}
D^{(l+1)}_{\{1\}}(u_{l+1})\cdots D^{(k)}_{\{1\}}(u_k)  \widehat\psi_{k}(u_1,\dots,u_k)(X)
(v_1\otimes v_2)
\\
\qquad {} =
\widehat\psi_{l}(u_1,\dots,u_k)\bigl(X D^{(l+1)}_{\{1\}}(u_{l+1})\cdots D^{(k)}_{\{1\}}(u_k)\bigr)
 (v_1 \otimes v_2)
\\
\qquad{} = \prod_{j=l+1}^k\big\langle T_{22}(u_j) v_1\big\rangle
\widehat\psi_{l}(u_1,\dots,u_k)(X)(v_1 \otimes v_2) .
\end{gather*}

Recall that we regard $\widehat\psi_{l}(u_1,\dots,u_k)(X)$ as a matrix over $ Y({\mathfrak{gl}}_N)\otimes Y({\mathfrak{gl}}_N)$.
All entries of this matrix belong to
$\psi\bigl( Y({\mathfrak{gl}}_{N-1}))\bigr)\otimes\psi\bigl( Y({\mathfrak{gl}}_{N-1}))\bigr)$. It follows from
relations~\eqref{Tabcd} that for any $X'\in Y({\mathfrak{gl}}_{N-1})$ the coef\/f\/icients of
the commutator $T_{11}(u)\psi(X')-\psi(X')T_{11}(u)$ belong to
the left ideal in $ Y({\mathfrak{gl}}_N)$ generated by the coef\/f\/icients of the series $T_{21}(u),\dots, T_{N1}(u)$.
Therefore,
\[
A(u_i) \psi(X') v_2 =\psi(X') T_{11}(u_i) v_2=
\big\langle T_{11}(u_i) v_2\big\rangle \psi(X')v_2
\]
because $A(u_i)=T_{11}(u_i)$, cf.~\eqref{ABD}, and $v_2$ is a weight singular
vector. Hence,
\begin{gather*}
A_{\{2\}}(u_1)\cdots A_{\{2\}}(u_l) \widehat\psi_{l}(u_1,\dots,u_k)(X)(v_1\otimes v_2)
\\
\qquad {}= \prod_{i=1}^l \big\langle T_{11}(u_i) v_2\big\rangle
\widehat\psi_{l}(u_1,\dots,u_k)(X)(v_1\otimes v_2),
\end{gather*}
which proves the lemma.
\end{proof}

\begin{prop}
\label{Step}
In the notation of Theorem~\ref{BBvox}, we have
\begin{gather}
 {\mathbb B}_{\xi}(t)(v_1\otimes v_2)
\label{step}
\\
\nonumber
{} = \sum_{l=0}^{\xi^1}\frac1{l!(\xi^1-l)!}
\operatorname{Sym}_{t^1_1,\dots, t^1_{\xi^1}}^{(1,\dots,\xi^1)}
\Biggl( \prod_{1\leqslant i<j\leqslant\xi^1} \frac{t^1_i-t^1_j-1}{t^1_i-t^1_j}
\prod_{i=1}^l \big\langle T_{11}\big(t^1_i\big)v_2\big\rangle
\prod_{j=l+1}^{\xi^1}\big\langle T_{22}\big(t^1_j\big)v_1\big\rangle
\\
\nonumber
{}\times B^{[1]}_{\{1\}}\big(t^1_1\big)\cdots B^{[l]}_{\{1\}}\big(t^1_l\big)
B^{[l+1]}_{\{2\}}\big(t^1_{l+1}\big)\cdots B^{[\xi^1]}_{\{2\}}\big(t^1_{\xi^1}\big) \widehat\psi_{l}\big(t^1_1,\dots, t^1_{\xi^1}\big)
\bigl({\mathbb B}^{\langle N-1\rangle}_{\ddot\xi}(\ddot t)\bigr)\Biggr)(v_1\otimes v_2).
\end{gather}
Here the space ${V_1\otimes V_2}$ is considered as the $Y({\mathfrak{gl}}_N)$-module in the left side of
the formula, and as the $Y({\mathfrak{gl}}_N)\otimes Y({\mathfrak{gl}}_N)$-module in the right side.
\end{prop}

\begin{proof}
Expand ${\mathbb B}_{\xi}(t)$ according to formula \eqref{BB1i}.
Since $ Y({\mathfrak{gl}}_N)$ acts in $V_1\otimes V_2$ via the copro\-duct~$\Delta$, we have
\begin{gather*}
{\mathbb B}_{\xi}(t)(v_1\otimes v_2)
 = \Delta\Bigl(B^{[1]}(t_1)\cdots B^{[\xi^1]}(t_{\xi^1}) \psi\big(t^1_1,\dots, t^1_{\xi^1}\big)
\bigl({\mathbb B}^{\langle N-1\rangle}_{\ddot\xi}(\ddot t)\bigr)\Bigr)(v_1 \otimes v_2)
\\
\hphantom{{\mathbb B}_{\xi}(t)(v_1\otimes v_2)}{}
=\Delta\bigl(B^{[1]}(t_1)\cdots B^{[\xi^1]}(t_{\xi^1})\bigr) \Delta\Bigl(
\psi\big(t^1_1,\dots, t^1_{\xi^1}\big)\bigl({\mathbb B}^{\langle N-1\rangle}_{\ddot\xi}(\ddot t)\bigr)\Bigr)(v_1\otimes v_2).
\end{gather*}
Recall that $\Delta$ applies to matrices entrywise. In the last expression, we
develop the fac\-tor $\Delta\bigl(B^{[1]}(t_1)\cdots B^{[\xi^1]}(t_{\xi^1})\bigr)$
according to Proposition~\ref{DlB}, and replace the factor
$\Delta\bigl(\psi\big(t^1_1,\dots$,\linebreak  $t^1_{\xi^1}\big)\bigl({\mathbb B}^{\langle N-1\rangle}_{\ddot\xi}(\ddot t)\bigr)\bigr)$
by $\widehat\psi_{\xi^1}\big(t^1_1,\dots, t^1_{\xi^1}\big)\bigl({\mathbb B}^{\langle N-1\rangle}_{\ddot\xi}(\ddot t)\bigr)$
according to Lemma~\ref{DlN}. After that, we utilize Lemma~\ref{ADX}
to transform the result to the right side of formula~\eqref{step}.
\end{proof}

\section{Trigonometric weight functions}\label{:Twtf}
Notation in this section may not coincide with the notation
in Sections~\ref{:Rwtf}--\ref{:SBBvox}.

The quantum loop algebra $U_q(\widetilde{{\mathfrak{gl}}_N})$ (the quantum af\/f\/ine algebra without central extension) is the unital
associative algebra with generators $L_{ab}^{\{\pm s\}}$, $a,b=1,\dots,N$, and
$s=0,1,2,\ldots$. Organize them into generating series
\begin{gather}
L^{\pm}_{ab}(u)=L_{ab}^{\{\pm 0\}}+
\sum_{s=1}^\infty L_{ab}^{\{\pm s\}} u^{\pm s},\qquad a,b=1,\dots,N,
\label{Lab}
\end{gather}
and combine the series into matrices
$L^{\pm}(u) =\sum\limits_{a,b=1}^{N}E_{ab}\otimes L^{\pm}_{ab}(u)$.
The def\/ining relations in $U_q(\widetilde{{\mathfrak{gl}}_N})$ are
\begin{gather*}
L_{ab}^{\{+0\}}= L_{ba}^{\{-0\}}=0, \qquad 1\leqslant a<b\leqslant N,
\\
L_{aa}^{\{-0\}} L_{aa}^{\{+0\}}=L_{aa}^{\{+0\}}L_{aa}^{\{-0\}}=1,
\qquad a=1,\dots,N,
\\
R_q^{(12)}(u/v)\bigl(L^\mu(u)\bigr)^{(1)}\bigl(L^\nu(v)\bigr)^{(2)}
= \bigl(L^\nu(v)\bigr)^{(2)}\bigl(L^\mu(u)\bigr)^{(1)}R_q^{(12)}(u/v),
\end{gather*}
$(\mu,\nu)=({+},{+}),({+},{-}),({-},{-})$,
where $R_q(u)$ is the trigonometric $R$-matrix~\eqref{Rq}.

The quantum loop algebra $U_q(\widetilde{{\mathfrak{gl}}_N})$ is a Hopf algebra. In terms of generating series
\eqref{Lab}, the coproduct $\Delta:U_q(\widetilde{{\mathfrak{gl}}_N})\to U_q(\widetilde{{\mathfrak{gl}}_N})\otimes U_q(\widetilde{{\mathfrak{gl}}_N})$ reads as follows:
\[
\Delta: \ L^{\pm}_{ab}(u)\mapsto \sum_{c=1}^N L^{\pm}_{cb}(u)\otimes L^{\pm}_{ac}(u).
\]
The subalgebras $U^\pm_q(\widetilde{{\mathfrak{gl}}_N})\subset U_q(\widetilde{{\mathfrak{gl}}_N})$ generated by the coef\/f\/icients of
the respective series $L^{\pm}_{ab}(u)$, $a,b=1,\dots,N$, are Hopf subalgebras.

There is a one-parameter family of automorphisms $\rho_x:U_q(\widetilde{{\mathfrak{gl}}_N})\to U_q(\widetilde{{\mathfrak{gl}}_N})$,
def\/ined by the rule
\[
\rho_x: \ L^{\pm}_{ab}(u)\mapsto L^{\pm}_{ab}(u/x).
\]

The quantum loop algebra $U_q(\widetilde{{\mathfrak{gl}}_N})$ contains the algebra $U_q({\mathfrak{gl}}_N)$ as a Hopf subalgebra.
The subalgebra is generated by the elements
$L^{\{+0\}}_{ab}$, $L^{\{-0\}}_{ab}$, $1\leqslant a\leqslant b\leqslant N$.
Set $\hat k_a=L^{\{-0\}}_{aa}$, $a=1,\dots,N$, and
\begin{gather}
\hat e_{ab}=-\frac{L^{\{+0\}}_{ba}\hat k_a}{q-q^{-1}},\qquad
\hat e_{ba}=\frac{\hat k_a^{-1}L^{\{-0\}}_{ab}}{q-q^{-1}},\qquad 1\leqslant a<b\leqslant N.
\label{eh}
\end{gather}
The elements $\hat k_1,\dots,\hat k_N$, $\hat e_{12},\dots,\hat e_{N-1,N}$,
$\hat e_{21},\dots,\hat e_{N,N-1}$ are the Chevalley generators \linebreak of~$U_q({\mathfrak{gl}}_N)$.
We list some of relations for the introduced elements below, subscripts
running over all possible values unless the range is specif\/ied explicitly:
\begin{gather*}
 \hat k_a \hat e_{bc}= q^{\delta_{ab}-\delta_{ac}}\hat e_{bc}\hat k_a,
\\
\begin{array}{@{}l}
 \hat e_{a,a+1}\hat e_{a+1,b}-q\hat e_{a+1,b}\hat e_{a,a+1}=
\hat e_{a,b-1}\hat e_{b-1,b}-q\hat e_{b-1,b}\hat e_{a,b-1}=\hat e_{ab},
\\[0.5ex]
 \hat e_{b,a+1}\hat e_{a+1,a}-q^{-1}\hat e_{a+1,a}\hat e_{b,a+1}=
\hat e_{b,b-1}\hat e_{b-1,a}-q^{-1}\hat e_{b-1,a}\hat e_{b,b-1}=\hat e_{ba},
\end{array}
\qquad a+1<b,
\\
 \hat e_{ca}\hat e_{ba}=q\hat e_{ba}\hat e_{ca},\qquad
\hat e_{cb}\hat e_{ca}=q\hat e_{ca}\hat e_{cb}, \qquad a<b<c.
\end{gather*}
The coproduct formulae are $\Delta(\hat k_a)=\hat k_a\otimes\hat k_a$,
\begin{gather*}
\Delta(\hat e_{a,a+1})=1\otimes \hat e_{a,a+1}+\hat e_{a,a+1}\otimes \hat k_a \hat k_{a+1}^{-1},
\qquad
\Delta(\hat e_{a+1,a})=\hat e_{a+1,a}\otimes 1
+\hat k_{a+1}\hat k_a^{-1}\otimes \hat e_{a+1,a}.
\end{gather*}

By minor abuse of notation we say that a vector $v$ in a $U_q({\mathfrak{gl}}_N)$-module has weight
$\big(\Lambda^{1},\dots$, $\Lambda^{N}\big)$ if $\hat k_a v=q^{\Lambda^{a}} v$ for all $a=1,\dots,N$. A vector~$v$
is called a \textit{singular vector} if $\hat e_{ba}v=0$ for all ${1\leqslant a<b\leqslant N}$.

The evaluation homomorphism $\epsilon: U_q(\widetilde{{\mathfrak{gl}}_N})\to U_q({\mathfrak{gl}}_N)$ is given by the rule
\begin{alignat*}{4}
&   L^+_{aa}(u) \mapsto \hat k_a^{-1} -u \hat k_a , \qquad &&
 L^-_{aa}(u) \mapsto \hat k_a-u^{-1}\hat k_a^{-1},\qquad && a=1,\dots,N, &
\\
&  L^+_{ab}(u)\mapsto -u\big(q-q^{-1}\big)\hat k_a \hat e_{ba}, \qquad &&
 L^-_{ab}(u)\mapsto \big(q-q^{-1}\big) \hat k_a \hat e_{ba},  \qquad &&&
\\
&  L^+_{ba}(u) \mapsto -\big(q-q^{-1}\big)\hat e_{ab}\hat k_a^{-1},\qquad &&
 L^-_{ba}(u)\mapsto u^{-1}\big(q-q^{-1}\big) \hat e_{ab}\hat k_a^{-1}.\qquad&& \mbox{\raisebox{12pt}[0pt][0pt]{$1\leqslant a<b\leqslant N,$}}
\end{alignat*}
Both the automorphisms $\rho_x$ and the homomorphism $\epsilon$ restricted to the subalgebra
$U_q({\mathfrak{gl}}_N)$ are the identity maps.

For a $U_q({\mathfrak{gl}}_N)$-module $V$ denote by $V(x)$ the $U_q(\widetilde{{\mathfrak{gl}}_N})$-module induced from $V$ by the homomorphism
${\epsilon\circ \rho_x}$. The module $V(x)$ is called an \textit{evaluation module} over $U_q(\widetilde{{\mathfrak{gl}}_N})$.
\begin{rem}
In a $k$-fold tensor product of evaluation modules the series $ L^+(u)$ and $ L^-(u)$
act as polynomials in $u$ and $u^{-1}$, respectively, and the action of $ L^+(u)$ is
proportional to that of $u^k L^-(u)$.
\end{rem}

Let $V$ be a $U^-_q(\widetilde{{\mathfrak{gl}}_N})$-module. A vector ${v\in V}$ is called a \textit{weight singular
vector} with respect to the action of $U^-_q(\widetilde{{\mathfrak{gl}}_N})$ if $ L^-_{ba}(u)v=0$ for all $1\leqslant a<b\leqslant N$,
and $v$ is an eigenvector for the action of $ L^-_{11}(u),\dots,L^-_{NN}(u)$;
the respective eigenvalues are denoted by
$\langle L^-_{11}(u)v\rangle,\dots, \langle L^-_{NN}(u) v\rangle$.

\begin{Example}
Let $V$ be a $U_q({\mathfrak{gl}}_N)$-module and let ${v\in V}$ be a singular vector of weight
$\big(\Lambda^{1},\dots,\Lambda^{N}\big)$. Then $v$ is a weight singular vector with respect to the action of $U^-_q(\widetilde{{\mathfrak{gl}}_N})$
in the evaluation modu\-le~$V(x)$ and $\langle L^-_{aa}(u)v\rangle=
q^{\Lambda^{a}}-q^{-\Lambda^{a}}x u^{-1}$, $a=1,\dots,N$.
\end{Example}

We will use two embeddings of the algebra $U_q(\widetilde{{\mathfrak{gl}}_{N-1}})$ into $U_q(\widetilde{{\mathfrak{gl}}_N})$,
called $\phi$ and $\psi$:
\begin{gather}
\phi \Bigl(\bigl(L^{\pm}_{ab}(u)\bigr)^{\langle N-1\rangle}\Bigr)=
\bigl(L^{\pm}_{ab}(u)\bigr)^{\langle N\rangle}, \qquad
\psi\Bigl(\bigl(L^{\pm}_{ab}(u)\bigr)^{\langle N-1\rangle}\Bigr)=
\bigl(L^{\pm}_{a+1,b+1}(u)\bigr)^{\langle N\rangle}.
\label{embeq}
\end{gather}
Here $\bigl(L^{\pm}_{ab}(u)\bigr)^{\langle N-1\rangle}$ and $\bigl(L^{\pm}_{ab}(u)\bigr)^{\langle N\rangle}$
are series~\eqref{Lab} for the algebras $U_q(\widetilde{{\mathfrak{gl}}_{N-1}})$ and $U_q(\widetilde{{\mathfrak{gl}}_N})$, respectively.

The constructions and statements in the rest of the section are similar to
those of Section~\ref{:Rwtf}. We will mention only essential points and omit
details.

Let $k$ be a nonnegative integer.
Let ${\xi=(\xi^1,\dots,\xi^{N-1})}$ be a collection of nonnegative integers.
Remind that ${\xi^{<a}=\xi^1+\cdots+\xi^{a-1}}$, $a=1,\dots,N$,
and ${|\xi|=\xi^1+\cdots +\xi^{N-1}=\xi^{<N}}$.

Consider a series in $|\xi|$ variables $t^1_1,\dots, t^1_{\xi^1},\dots, t^{N-1}_1,\dots,t^{N-1}_{\xi^{N-1}}$ with coef\/f\/icients in
$U^-_q(\widetilde{{\mathfrak{gl}}_N})$:
\begin{gather}
{\widehat{\mathbb B}}_{\xi}\big(t^1_1,\dots, t^{N-1}_{\xi^{N-1}}\big)=
\big({\operatorname{tr}^{\otimes|\xi|}} \otimes {\rm id}\big)\Biggl(
\bigl( L^-\big(t^1_1\big)\bigr)^{(1,|\xi|+1)}\cdots
\bigl( L^-\big(t^{N-1}_{\xi^{N-1}}\big)\bigr)^{(|\xi|,|\xi|+1)}
\nonumber\\
\hphantom{{\widehat{\mathbb B}}_{\xi}\big(t^1_1,\dots, t^{N-1}_{\xi^{N-1}}\big)=}{}\times
\mathop{\overrightarrow\prod}\limits_{(a,i)<(b,j)}
R_q^{(\xi^{<b}+j,\xi^{<a}+i)}\big(t^b_j/t^{a}_i\big)
E_{21}^{\otimes\xi^1} \otimes\cdots\otimes E_{N,N-1}^{\otimes\xi^{N-1}}\otimes 1\Biggr),\label{BBqh}
\end{gather}
the convention being the same as in \eqref{BBh}.

\begin{rem}
The series ${\widehat{\mathbb B}}_{\xi}\big(t^1_1,{\dots}, t^{N-1}_{\xi^{N-1}}\big)$ belongs to $U_q(\widetilde{{\mathfrak{gl}}_N})\bigl[t^1_1,{\dots}, t^{N-1}_{\xi^{N-1}}\bigr]\bigl[ \bigl[\big(t^1_1\big)^{-1}\!,{\dots},\big(t^{N-1}_{\xi^{N-1}}\big)^{-1}\bigr]\bigr]$.
\end{rem}

Set
\begin{gather}
{\mathbb B}_{\xi}(t)={\widehat{\mathbb B}}_{\xi}(t)\prod_{a=1}^{N-1}
\prod_{1\leqslant i<j\leqslant\xi^a}\frac{t^{a}_i}{qt^{a}_j-q^{-1}t^{a}_i}
\prod_{1\leqslant a<b<N}\prod_{i=1}^{\xi^a}\prod_{j=1}^{\xi^b}
\frac{t^{a}_i}{t^b_j-t^{a}_i},
\label{BBq}
\end{gather}
To indicate the dependence on $N$, if necessary,
we will write ${\mathbb B}^{\langle N\rangle}_{\xi}(t)$.

\begin{Example}
Let $N=2$ and $\xi=(\xi^1)$. Then
${\mathbb B}^{\langle 2\rangle}_{\xi}(t)= L^-_{12}\big(t^1_1\big)\cdots L^-_{12}\big(t^1_{\xi^1}\big)$.
\end{Example}

\begin{Example}
Let $N=3$ and $\xi=(1,1)$. Then
\[
{\mathbb B}^{\langle 3\rangle}_{\xi}(t)= L^-_{12}\big(t^1_1\big) L^-_{23}\big(t^2_1\big)+
\big(q-q^{-1}\big) \frac{t^2_1}{t^2_1-t^1_1} L^-_{13}\big(t^1_1\big) L^-_{22}\big(t^2_1\big).
\]
\end{Example}

\begin{Example}
Let $N=4$ and $\xi=(1,1,1)$. Then
\begin{gather*}
{\mathbb B}^{\langle 4\rangle}_{\xi} (t)=
 L^-_{12}\big(t^1_1\big) L^-_{23}\big(t^2_1\big) L^-_{34}\big(t^3_1\big)
\\
\hphantom{{\mathbb B}^{\langle 4\rangle}_{\xi} (t)=}{}
+\big(q-q^{-1}\big)\left(\frac{t^2_1}{t^2_1-t^1_1}
 L^-_{13}\big(t^1_1\big) L^-_{22}\big(t^2_1\big) L^-_{34}\big(t^3_1\big)
 + \frac{t^3_1}{t^3_1-t^2_1} L^-_{12}\big(t^1_1\big)
 L^-_{24}\big(t^2_1\big)  L^-_{33}\big(t^3_1\big) \right)
\\
\hphantom{{\mathbb B}^{\langle 4\rangle}_{\xi} (t)=}{}
+\big(q-q^{-1}\big)^2\frac{t^2_1t^3_1}{\big(t^2_1-t^1_1\big)\big(t^3_1-t^2_1\big)}
\bigl(L^-_{14}\big(t^1_1\big) L^-_{22}\big(t^2_1\big) L^-_{33}\big(t^3_1\big)+
 L^-_{13}\big(t^1_1\big) L^-_{24}\big(t^2_1\big) L^-_{32}\big(t^3_1\big)\bigr)
\\
\hphantom{{\mathbb B}^{\langle 4\rangle}_{\xi} (t)=}{}
+\big(q-q^{-1}\big)t^3_1
\frac{\big(t^2_1-t^1_1\big)\big(t^3_1-t^2_1\big)+\big(q-q^{-1}\big)^2 t^2_1t^3_1}
{\big(t^2_1-t^1_1\big) \big(t^3_1-t^1_1\big)\big(t^3_1-t^2_1\big)}
 L^-_{14}\big(t^1_1\big) L^-_{23}\big(t^2_1\big) L^-_{32}\big(t^3_1\big).
\end{gather*}
\end{Example}

Recall that the direct product of the symmetric groups $S_{\xi^1}\times \cdots \times S_{\xi^{N-1}}$ acts on expressions
in $|\xi|$ variables, permuting the variables with the same superscript,
cf.~\eqref{sixi}.

\begin{lem}[\protect{\cite[Theorem~3.3.4]{TV1}}]
\label{BBSq}
The expression ${\mathbb B}_{\xi}(t)$ is invariant under the action of the group
${S_{\xi^1}\times \cdots \times S_{\xi^{N-1}}}$.
\end{lem}

If $v$ is a weight singular vector with respect to the action of $U^-_q(\widetilde{{\mathfrak{gl}}_N})$, we call the
expression ${\mathbb B}_{\xi}(t)v$ a (\textit{trigonometric}) \textit{vector-valued weight function} of weight $\big(\xi^1,{\xi^2-\xi^1},\dots,{\xi^{N-1}-\xi^{N-2}}$, $-\xi^{N-1}\big)$ associated with $v$.

Weight functions associated with $U_q({\mathfrak{gl}}_N)$ weight singular vectors in evaluation $U^-_q(\widetilde{{\mathfrak{gl}}_N})$-mod\-ules
(in particular, highest weight vectors of highest weight $U_q({\mathfrak{gl}}_N)$-modules) can be calculated explicitly
by means of the following Theorems~\ref{BBNq} and~\ref{BB1q}, which are analogues
of Theorems~\ref{BBN} and \ref{BB1}, respectively. Corollaries~\ref{BCNq} and~\ref{BC1q} are
the respective counterparts of Corollaries~\ref{BCN} and~\ref{BC1}.

Theorem~\ref{BBvoxq} and Corollary~\ref{BBvoxnq} are analogous to Theorem~\ref{BBvox}
and Corollary~\ref{BBvoxn} in the Yangian case and yield combinatorial formulae
for weight functions associated with tensor products of highest weight vectors of highest weight evaluation modules
over the quantum loop algebra.

\begin{rem}
The expression for a vector-valued weight function used here may dif\/fer from
the expressions for the corresponding objects used in other papers,
see~\cite{KR, TV1}. The discrepancy can occur due to the choice
of coproduct for the quantum loop algebra $U_q(\widetilde{{\mathfrak{gl}}_N})$ as well as the choice of normalization.
\end{rem}

For a nonnegative integer $k$ introduce a function $W_{k}(t_1,\dots,t_k)$:
\[
W_{k}(t_1,\dots,t_k)=
\prod_{1\leqslant i<j\leqslant k}\frac{q^{-1}t_i-qt_j}{t_i-t_j}.
\]
For an expression $f(t^1_1,\dots, t^{N-1}_{\xi^{N-1}})$, set
\begin{gather}
\overline{\operatorname{Sym}}^\xi_{t}f(t)=\operatorname{Sym}^\xi_{t}\left(f(t)
\prod_{a=1}^{N-1}W_{\xi^a}\big(t^{a}_1,\dots,t^{a}_{\xi^a}\big)\right),
\label{Symbq}
\end{gather}
where $\operatorname{Sym}^\xi_{t}$ is def\/ined by \eqref{Sym}.

Let ${\eta^1\leqslant\cdots\leqslant\eta^{N-1}}$ be nonnegative integers.
Def\/ine a function $X_\eta\big(t^1_1,\dots, t^1_{\eta^1};\dots; t^{N-1}_1$, $\dots, t^{N-1}_{\eta^{N-1}}\big)$,
\begin{gather}
X_\eta(t)=\prod_{a=1}^{N-2}\left[\prod_{j=1}^{\eta^a}
\frac1{t^{a+1}_j-t^{a}_j}\prod_{i=1}^{j-1}
\frac{qt^{a+1}_i-q^{-1}t^{a}_j}{t^{a+1}_i-t^{a}_j}\right].
\label{Xq}
\end{gather}
The function $X_\eta(t)$ does not actually depend on the variables
$t^{N-1}_{\eta^{N-2}+1},\dots,t^{N-1}_{\eta^{N-1}}$.

For nonnegative integers ${\eta^{1}\geqslant\cdots \geqslant\eta^{N-1}}$
def\/ine a function $Y_{\eta}\big(t^1_1,\dots, t^1_{\eta^1};\dots; t^{N-1}_1,\dots$,\linebreak  $t^{N-1}_{\eta^{N-1}}\big)$,
\begin{gather}
Y_{\eta}(t)=\prod_{a=2}^{N-1}\left[\prod_{j=1}^{\eta^a}
\frac1{t^{a}_j-t^{a-1}_{j+\eta^{a-1}-\eta^a}}
\prod_{i=1}^{j-1}
\frac{q t^{a}_i-q^{-1}t^{a-1}_{j+\eta^{a-1}-\eta^a}}
{t^{a}_i-t^{a-1}_{j+\eta^{a-1}-\eta^a}}\right].
\label{Yq}
\end{gather}
The function $Y_{\eta}(t)$ does not actually depend on the variables
$t^1_1,\dots,t^1_{\eta^1-\eta^2}$.

For any $\xi,\eta\in{\mathbb Z}_{\geqslant 0}^{N-1}$,
def\/ine a function $Z_{\xi,\eta}\big(t^1_1,\dots, t^{N-1}_{\xi^{N-1}};s^1_1,\dots,s^{N-1}_{\eta^{N-1}}\big)$,
\begin{gather}
Z_{\xi,\eta}(t;s)=\prod_{a=1}^{N-2}\prod_{i=1}^{\xi^{a+1}}
\prod_{j=1}^{\eta^a}\frac{qt^{a+1}_i-q^{-1}s^a_j}{t^{a+1}_i-s^a_j}.
\label{Zq}
\end{gather}
The function $Z_{\xi,\eta}(t;s)$ does not depend on the variables $t^1_1,\dots, t^1_{\xi^1}$ and
$s^{N-1}_1,\dots, s^{N-1}_{\eta^{N-1}}$.

We are using the following $q$-numbers:
$[n]_q=\frac{q^n-q^{-n}}{q-q^{-1}}$, and $q$-factorials:
$[n]_q!=\prod\limits_{r=1}^n[r]_q$.
Recall that for a collection $t$ of $|\xi|$ variables
we introduced the subcollections $t_{[\eta]}$, $t_{(\eta,\xi]}$ by~(\ref{dot})
and~$\dot t$,~$\ddot t$ by~(\ref{teta}).

For any $1\leqslant a<b\leqslant N$ set $\check e_{ba}=\hat k_a \hat e_{ab}$, cf.~\eqref{eh}.

\begin{thm}
\label{BBNq}
Let $V$ be a $U_q({\mathfrak{gl}}_N)$-module and ${v\in V}$ a singular vector of weight $\big(\Lambda^{1},\dots,\Lambda^{N}\big)$.
Let $\xi^1,\dots,\xi^{N-1}$ be nonnegative integers and $t=\big(t^1_1,\dots, t^1_{\xi^1};\dots; t^{N-1}_1,\dots, t^{N-1}_{\xi^{N-1}}\big)$.
In the eva\-luation $U_q(\widetilde{{\mathfrak{gl}}_N})$-module $V(x)$, one has
\begin{gather}
 {\mathbb B}_{\xi}(t) v =
\sum_\eta\big(q-q^{-1}\big)^{|\eta|} \frac1{[\eta^{1}]_q!}
\prod_{a=1}^{N-2}\frac{q^{\eta^a (\eta^a -\eta^{a+1})}}
{[\xi^a -\eta^{a}]_q![\eta^{a+1} -\eta^{a}]_q!}\nonumber\\
\hphantom{{\mathbb B}_{\xi}(t) v =}{}
\times \overline{\operatorname{Sym}}^\xi_{t}
\Biggl[ X_\eta(t_{(\xi-\eta,\xi]})
Z_{\xi-\eta,\eta}(t_{[\xi-\eta]};t_{(\xi-\eta,\xi]})
\prod_{a=1}^{N-2}\prod_{i=0}^{\eta^a-1}\bigl(
q^{\Lambda^{ a+1}} t^{a}_{\xi^a -i} -q^{-\Lambda^{a+1}}x\bigr)\nonumber\\
\hphantom{{\mathbb B}_{\xi}(t) v =}{}
\times \check e_{N,N-1}^{ \eta^{N-1} -\eta^{N-2}}
\check e_{N,N-2}^{\eta^{N-2} -\eta^{N-3}} \cdots\check e_{N1}^{ \eta^1}
\phi\bigl({\mathbb B}^{\langle N-1\rangle}_{(\xi-\eta)^{\cdot}}(\dot t_{[\xi-\eta]})\bigr) v
\Biggr],\label{BBXq}
\end{gather}
the sum being taken over all $\eta=\big(\eta^{1},\dots,\eta^{N-1}\big)\in{\mathbb Z}_{\geqslant 0}^{N-1}$ such that
$\eta^{1}\leqslant\cdots \leqslant\eta^{N-1}=\xi^{N-1} $ and $\eta^{a} \leqslant\xi^a$
for all $a=1,\dots,N-2$. Other notation is as follows: $\overline{\operatorname{Sym}}^\xi_{t}$
is defined by \eqref{Symbq}, the functions $X_\eta$ and
$Z_{\xi-\eta,\eta}$ are respectively given by formulae \eqref{Xq} and \eqref{Zq},
$\phi$ is the first of embeddings \eqref{embeq}, and
\[
{{\mathbb B}^{\langle N-1\rangle}_{(\xi-\eta)^{\cdot}(\dot t_{[\xi-\eta]})} =
{\mathbb B}^{\langle N-1\rangle}_\zeta(s)\big|_{\zeta=(\xi-\eta)^{\cdot},\; s=\dot t_{[\xi-\eta]}}},
\]
${\mathbb B}^{\langle N-1\rangle}_\zeta(s)$ coming from \eqref{BBq}.
\end{thm}

\begin{rem}
For ${N=2}$, the sum in the right side of formula \eqref{BBXq} contains only one term:
${\eta=\xi}$. Moreover, $X_\eta=Z_{\xi-\eta,\eta}=1$, and
${\mathbb B}^{\langle1\rangle:}_{(\xi-\eta)^{\cdot}}=1$ by convention.
\end{rem}

\begin{cor}
\label{BCNq}
Let $V$ be a $U_q({\mathfrak{gl}}_N)$-module and ${v\in V}$ a singular vector of weight $\big(\Lambda^{1},\dots,\Lambda^{N}\big)$.
Let $\xi^1,\dots,\xi^{N-1}$ be nonnegative integers and $t=(t^1_1,\dots, t^1_{\xi^1};\dots; t^{N-1}_1,\dots, t^{N-1}_{\xi^{N-1}})$.
In the eva\-luation $U_q(\widetilde{{\mathfrak{gl}}_N})$-module $V(x)$, one has
\begin{gather}
{\mathbb B}_{\xi}(t) v = \big(q-q^{-1}\big)^{|\xi|}
 \sum_m \left[\mathop{\overleftarrow\prod}\limits_{1\leqslant b<a\leqslant N}
\frac{q^{m^{ a,b-1} (m^{ a,b-1} -m^{ ab})}}{[m^{ab} -m^{a,b-1}]_q!}
\check e_{ab}^{ m^{ ab} -m^{ a,b-1}} \right]v
\label{BcNq}\\
\hphantom{{\mathbb B}_{\xi}(t) v =}{}
\times
\overline{\operatorname{Sym}}^\xi_{t} \left[
\prod_{a=3}^N\prod_{b=1}^{a-2} \prod_{i=1}^{m^{ ab}}\left(
\frac{q^{\Lambda^{b+1}} t^b_{i+{\widetilde m}^{ab}} -q^{-\Lambda^{b+1}}x}
{t^{b+1}_{i+{\widetilde m}^{a,b+1}} -t^b_{i+{\widetilde m}^{ ab}}}
\!\prod_{1\leqslant j<i+{\widetilde m}^{ a,b+1}}\!
\frac{qt^{b+1}_j -q^{-1}t^b_{i+{\widetilde m}^{ ab}}}
{t^{b+1}_j -t^b_{i+{\widetilde m}^{ ab}}}\right)\right].\nonumber
\end{gather}
Here the sum is taken over all collections of nonnegative integers $m^{ab}$,
$1\leqslant b<a\leqslant N$, such that $m^{a1}\leqslant\cdots\leqslant m^{a,a-1}$ and
$m^{a+1,a}+\cdots +m^{Na}=\xi^a$ for all $a=1,\dots,N-1$; by convention,
$m^{a0}=0$ for any $a=2,\dots,N$. Other notation is as follows:
in the ordered product the fac\-tor~$\check e_{ab}^\circledast$ is to the left of the factor
$\check e_{cd}^\circledast$ if $a>c$, or $a=c$ and $b>d$, $\overline{\operatorname{Sym}}^\xi_{t}$ is defined
by~\eqref{Symbq}, and ${\widetilde m}^{ab}=m^{b+1,b}+\cdots +m^{a-1,b}$ for all
$1\leqslant b<a\leqslant N$, in particular, ${\widetilde m}^{a,a-1}=0$.
\end{cor}

\begin{thm}
\label{BB1q}
Let $V$ be a $U_q({\mathfrak{gl}}_N)$-module and ${v\in V}$ a singular vector of weight $\big(\Lambda^{1},\dots,\Lambda^{N}\big)$.
Let $\xi^1,\dots,\xi^{N-1}$ be nonnegative integers and $t=\big(t^1_1,\dots, t^1_{\xi^1};\dots; t^{N-1}_1,\dots, t^{N-1}_{\xi^{N-1}}\big)$.
In the eva\-luation $U_q(\widetilde{{\mathfrak{gl}}_N})$-module $V(x)$, one has
\begin{gather}
{\mathbb B}_{\xi}(t) v =
\sum_\eta \big(q-q^{-1}\big)^{|\eta|}\frac{1}{[\eta^{N-1}]_q!}
\prod_{a=2}^{N-1}\frac{q^{\eta^a (\eta^{a-1} -\eta^a)}}
{[\xi^a -\eta^{a}]_q![\eta^{a-1} -\eta^{a}]_q!}\nonumber\\
\hphantom{{\mathbb B}_{\xi}(t) v =}{}
\times \overline{\operatorname{Sym}}^\xi_{t} \Biggl[Y_\eta(t_{[\eta]})
Z_{\eta,\xi-\eta}(t_{[\eta]};t_{(\eta,\xi]})
\prod_{a=2}^{N-1}\prod_{i=1}^{\eta^a}
\bigl(q^{\Lambda^{ a}} t^{a}_i -q^{-\Lambda^{ a}}x\bigr)\nonumber
\\
\hphantom{{\mathbb B}_{\xi}(t) v =}{}
\times \check e_{21}^{\eta^1-\eta^2}\check e_{31}^{\eta^2-\eta^3}
\cdots \check e_{N1}^{\eta^{N-1}}\psi\bigl({\mathbb B}^{\langle N-1\rangle}_{(\xi-\eta)^{\cdot\cdot}}
(\ddot t_{(\eta,\xi]})\bigr) v\Biggr], \label{BBYq}
\end{gather}
the sum being taken over all $\eta=\big(\eta^{1},\dots,\eta^{N-1}\big)\in{\mathbb Z}_{\geqslant 0}^{N-1}$ such that
$\xi^1 =\eta^{1}\geqslant \cdots\geqslant\eta^{N-1} $ and $\eta^{a}\leqslant\xi^a$
for all $a=2,\dots,N-1$. Other notation is as follows:
$\overline{\operatorname{Sym}}^\xi_{t}$ is defined by \eqref{Symbq}, the functions $Y_{\eta}$ and
$Z_{\eta,\xi-\eta}$ are respectively given by formulae~\eqref{Yq} and~\eqref{Zq},
$\psi$ is the second of embeddings~\eqref{embeq}, and
\[
{{\mathbb B}^{\langle N-1\rangle}_{(\xi-\eta)^{\cdot\cdot}}(\ddot t_{(\eta,\xi]}) =
{\mathbb B}^{\langle N-1\rangle}_\zeta(s)\big|_{\zeta=(\xi-\eta)^{\cdot\cdot},\; s=\ddot t_{(\eta,\xi]}}},
\]
${\mathbb B}^{\langle N-1\rangle}_\zeta(s)$ coming from \eqref{BBq}.
\end{thm}

\begin{rem}
For ${N=2}$, the sum in the right side of formula~\eqref{BBYq} contains only one term:
${\eta=\xi}$. Moreover, ${Y_\eta=Z_{\eta,\xi-\eta}=1}$, and
${\mathbb B}^{\langle 1\rangle}_{(\xi-\eta)^{\cdot\cdot}}=1$ by convention.
\end{rem}

\begin{cor}
\label{BC1q}
Let $V$ be a $U_q({\mathfrak{gl}}_N)$-module and ${v\in V}$ a singular vector of weight $\big(\Lambda^{1},\dots,\Lambda^{N}\big)$.
Let $\xi^1,\dots,\xi^{N-1}$ be nonnegative integers and $t=\big(t^1_1,\dots, t^1_{\xi^1};\dots; t^{N-1}_1,\dots, t^{N-1}_{\xi^{N-1}}\big)$.
In the eva\-luation $U_q(\widetilde{{\mathfrak{gl}}_N})$-module $V(x)$, one has
\begin{gather}
 {\mathbb B}_{\xi}(t) v =
\big(q-q^{-1}\big)^{|\xi|}\sum_m\left[\mathop{\overrightarrow\prod}\limits_{1\leqslant b<a\leqslant N}
\frac{q^{m^{ a+1,b} (m^{ ab} -m^{a+1,b})}}{[m^{ab} -m^{a+1,b}]_q!}
\check e_{ab}^{m^{ ab} -m^{a+1,b}} \right]v
\label{Bc1q}
\\
{} \times \overline{\operatorname{Sym}}^\xi_{t} \left[ \prod_{a=2}^{N-1}
\prod_{b=1}^{a-1}\prod_{i=0}^{m^{ a+1,b} -1}\left(
\frac{q^{\Lambda^{ a}} t^{a}_{{\widehat m}^{ a+1,b} -i} -q^{-\Lambda^{ a}}x}
{t^{a}_{{\widehat m}^{ a+1,b} -i} -t^{a-1}_{{\widehat m}^{ab} -i}}
\prod_{{\widehat m}^{ab} -i<j\leqslant\xi^{a-1}}
\frac{qt^{a}_{{\widehat m}^{a+1,b} -i} -q^{-1}t^{a-1}_j }
{t^{a}_{{\widehat m}^{ a+1,b} -i} -t^{a-1}_j}\right)\right].\nonumber
\end{gather}
Here the sum is taken over all collections of nonnegative integers $m^{ab}$,
$1\leqslant b<a\leqslant N$, such that $m^{a+1,a}\geqslant \cdots \geqslant m^{Na}$ and
$m^{a+1,1}+\cdots +m^{a+1,a} =\xi^a$ for all $a=1,\dots,N-1$; by convention,
$m^{N+1,a} =0$ for any $a=1,\dots,N$. Other notation is as follows: in
the ordered product the fac\-tor~$\check e_{ab}^\circledast$ is to the left of the fac\-tor~$\check e_{cd}^\circledast$ if $b<d$, or $b=d$ and $a<c$, $\overline{\operatorname{Sym}}^\xi_{t}$ is defined
by~\eqref{Symbq}, and ${\widehat m}^{ab} =m^{a1} +\cdots +m^{ab} $ for all
$1\leqslant b<a\leqslant N$, in particular, ${\widehat m}^{a+1,a} =\xi^a$.
\end{cor}

\begin{thm}[\cite{TV1}]
\label{BBvoxq}
Let  $V_1$, $V_2$ be $U^-_q(\widetilde{{\mathfrak{gl}}_N})$-modules and  ${v_1 \in V_1}$, ${v_2 \in V_2}$ weight
singular vectors with respect to the action of $U^-_q(\widetilde{{\mathfrak{gl}}_N})$. Let $\xi^1,\dots,\xi^{N-1}$ be nonnegative
integers and $t=\big(t^1_1,\dots, t^1_{\xi^1};\dots; t^{N-1}_1,\dots, t^{N-1}_{\xi^{N-1}}\big)$. Then
\begin{gather}
{\mathbb B}_{\xi}(t) (v_1  \otimes v_2)=
\sum_{\eta}\prod_{a=1}^{N-1}
\frac1{[\xi^a -\eta^{a}]_q![\eta^{a}]_q!}\overline{\operatorname{Sym}}^\xi_{t}
\Biggl[\prod_{a=1}^{N-2}\prod_{i=1}^{\eta^{a+1}}
\prod_{j=\eta^a +1}^{\xi^a}
\frac{qt^{a+1}_i-q^{-1}t^{a}_j}{t^{a+1}_i-t^{a}_j}\nonumber
\\
\times \prod_{a=1}^{N-1} \left(
\prod_{i=1}^{\eta^a}\big\langle  L^-_{a+1,a+1}\big(t^{a}_j\big) v_1\big\rangle
\prod_{j=\eta^a +1}^{\xi^a}
\big\langle  L^-_{aa}\big(t^{a}_i\big) v_2\big\rangle \right) {\mathbb B}_{\eta}(t_{[\eta]}) v_1
 \otimes{\mathbb B}_{\xi-\eta}(t_{(\eta,\xi]}) v_2\Biggr], \label{BBvxq}
\end{gather}
the sum being taken over all  $\eta=\big(\eta^{1},\dots,\eta^{N-1}\big)\in{\mathbb Z}_{\geqslant 0}^{N-1}$ such that
 $\xi-\eta\in{\mathbb Z}_{\geqslant 0}^{N-1}$. In the left side we assume that ${\mathbb B}_{\xi}(t)$
acts in the $U^-_q(\widetilde{{\mathfrak{gl}}_N})$-module $V_1\otimes V_2$.
\end{thm}

\begin{cor}
\label{BBvoxnq}
Let  $V_1,\dots, V_n$ be $U^-_q(\widetilde{{\mathfrak{gl}}_N})$-modules and  ${v_r \in V_r}$,  $r=1,\dots,n$, weight
singular vectors with respect to the action of $U^-_q(\widetilde{{\mathfrak{gl}}_N})$. Let $\xi^1,\dots,\xi^{N-1}$ be nonnegative
integers and $t=\big(t^1_1,\dots, t^1_{\xi^1};\dots; t^{N-1}_1,\dots, t^{N-1}_{\xi^{N-1}}\big)$. Then
\begin{gather}
 {\mathbb B}_{\xi}(t) ( v_1\otimes\cdots\otimes v_n)
\nonumber\\
\qquad {}  =\sum_{\eta_1,\dots, \eta_{n-1}} \prod_{a=1}^{N-1}
\prod_{r=1}^n \frac1{[\eta^{a}_r -\eta^{a}_{r-1}]_q!}
\overline{\operatorname{Sym}}^\xi_{t} \Biggl[\prod_{a=1}^{N-2}\prod_{r=1}^{n-1}
\prod_{i=\eta^{a+1}_{r-1}+1}^{\eta^{a+1}_r}
\prod_{j=\eta^a_r +1}^{\xi^a}
\frac{qt^{a+1}_i-q^{-1}t^{a}_j}{t^{a+1}_i-t^{a}_j}\nonumber
\\
\qquad\quad{}\times
\prod_{a=1}^{N-1}\prod_{r=1}^n\left(
\prod_{i=1}^{\eta^a_{r-1}}\big\langle L^-_{aa}\big(t^{a}_i\big) v_r\big\rangle
\prod_{j=\eta^a_r +1}^{\xi^a} \big\langle L^-_{a+1,a+1}\big(t^{a}_j\big) v_r\big\rangle
 \right) \nonumber
\\
\qquad\quad{}
\times {\mathbb B}_{\eta_1}(t_{[\eta_1 ]}) v_1 \otimes
{\mathbb B}_{\eta_2-\eta_1}(t_{(\eta_1 ,\eta_2]}) v_2 \otimes\cdots\otimes
{\mathbb B}_{\xi-\eta_{n-1}}(t_{(\eta_{n-1} ,\xi]}) v_n\Biggr].\label{BBvxnq}
\end{gather}
Here the sum is taken over all $\eta_1,\dots,\eta_{n-1}\in{\mathbb Z}_{\geqslant 0}^{N-1}$,
 $\eta_r=\big(\eta^{1}_r,\dots,\eta^{N-1}_r\big)$, such that
 $\eta_{r+1} -\eta_r\in{\mathbb Z}_{\geqslant 0}^{N-1}$ for any  $r=1,\dots,n-1$, and
 $\eta_0 =0$,  $\eta_n =\xi$,  by convention.
The sets $t_{[\eta_1 ]}$, $t_{(\eta_r ,\eta_{r+1}]}$ are defined by~\eqref{teta}. In the left side we assume that ${\mathbb B}_{\xi}(t)$ acts in the $U^-_q(\widetilde{{\mathfrak{gl}}_N})$-module
$ V_1\otimes\cdots\otimes V_n$.
\end{cor}

\begin{rem}
Denominators in the right sides of formulae \eqref{BBXq}--\eqref{BBvxnq} contain
$q$-factorials, which can vanish when $q$ is a root of unity. Nevertheless,
the right sides remain well-def\/ined at roots of unity. This happens due to the fact
that the symmetrized expressions in square brackets have nontrivial stationary
subgroups, cf.~Remark at the end of Section~\ref{:Rwtf}, so the result of
the symmetrization $\overline{\operatorname{Sym}}^\xi_{t}$ divided by the product of $q$-factorials
can be replaced by the sum over the cosets.
\end{rem}

Proofs of Theorems~\ref{BB1q}, \ref{BBNq} and \ref{BBvoxq} are similar to those
of Theorems~\ref{BB1}, \ref{BBN} and \ref{BBvox}, respectively. Here we mention only
the required modif\/ications of technical facts: identity \eqref{k!} and
Lemmas~\ref{Gyz},~\ref{Gyz1}. The analogue of the identity \eqref{k!} is
\begin{gather}
\sum_{\sigma \in S_k }\prod_{1\leqslant i<j\leqslant k}
\frac{q^{-1} s_{\sigma_i}-q s_{\sigma_j}}{s_{\sigma_i}-s_{\sigma_j}}=[k]_q!,
\label{k!q}
\end{gather}
and Lemmas~\ref{Gyz} and \ref{Gyz1} are to be generalized as follows.

\begin{lem}
\label{Gyzq}
Let $p$, $r$ be positive integers such that $p\leqslant r$. Then
\begin{gather}
 \operatorname{Sym}_{z_1,\dots,z_r}\left[
\prod_{i=1}^p \left( \frac1{y_i -z_{i+r-p}}\prod_{i<j\leqslant p}
\frac{qy_i -q^{-1} z_{j+r-p}}{y_i -z_{j+r-p}} \right)
\prod_{1\leqslant i<j\leqslant r} \frac{q^{-1} z_i -q z_j}{z_i -z_j}\right]
\nonumber
\\
=[r-p]_q!\sum_{{\boldsymbol{d}}}\operatorname{Sym}_{y_1,\dots,y_p}
\left[\prod_{i=1}^p \left( \frac{q^{i-d_i}}{y_i -z_{d_i}}
\prod_{d_i<j\leqslant r}\frac{qy_i -q^{-1} z_j}{y_i -z_j} \right)
\prod_{1\leqslant i<j\leqslant p}\! \frac{q^{-1}y_i -q y_j}{y_i -y_j}\right]\!,\!\!\!\label{yzq}
\\[1ex]
 \operatorname{Sym}_{z_1,\dots,z_r}\left[
\prod_{i=1}^p \left( \frac{q^{p-r}}{y_i -z_i}\prod_{1\leqslant j\leqslant i}
\frac{q^{-1}y_i -qz_j}{y_i -z_j} \right)
\prod_{1\leqslant i<j\leqslant r} \frac{q^{-1} z_i -q z_j}{z_i -z_j}\right]
\nonumber
\\
=[r-p]_q!\sum_{{\boldsymbol{d}}}\operatorname{Sym}_{y_1,\dots,y_p}
\left[\prod_{i=1}^p \left( \frac{q^{i-d_i}}{y_i -z_{d_i}}
\prod_{1\leqslant j<d_i}\frac{q^{-1}y_i -qz_j}{y_i -z_j} \right)
\prod_{1\leqslant i<j\leqslant p}\! \frac{q^{-1}y_i -q y_j}{y_i -y_j}\right]\!,\!\!\!\label{yzq1}
\end{gather}
the sums being taken over all $p$-tuples ${\boldsymbol{d}}=(d_1,\dots,d_p)$ such that
$1\leqslant d_1<\cdots <d_p\leqslant r$.
\end{lem}

\begin{proof}
It suf\/f\/ices to prove formula \eqref{yzq} since formulae \eqref{yzq} and \eqref{yzq1}
transform to each other by the change of variables $y_i \to y_{p-i}$,
$z_j \to z_{r-j}$, $q\to q^{-1}$, and a suitable change of summation
indices.

Consider the left side of  formula~\eqref{yzq} as a function of $z_1,\dots,z_r$ and denote it $f(z_1,\dots,z_r)$.
It has the following properties.
\begin{enumerate}\itemsep=0pt
\item[i)] $f(z_1,\dots,z_r)$ is symmetric in $z_1,\dots,z_r$.
\item[ii)]
$f(z_1,\dots,z_r)$ is a rational function of $z_1$ with only simple poles located at ${z_1 =y_i}$,
$i=1,\dots,p$, and regular as $z_1\to\infty$.
\item[iii)]
$\mathop{\operatorname{Res}}\limits_{z_1 =y_i} f\big(z_1,q^{2} y_i,z_3,\dots,z_r\big)=0$ for any $i=1,\dots,p$.
\item[iv)] $f(uz_1,\dots,uz_r)=u^{p-r}\bigl(1+o(1)\bigr)$ as $u\to\infty$.
\end{enumerate}

Denote by $C_{r-p}(y_1,\dots,y_p;z_1,\dots,z_r)$ the collection of properties  i)--iv),
the subscript $r-p$ referring to the exponent of $u$ in property~iv).

Consider the partial fraction decomposition of $f(z_1,\dots,z_r)$ as a function of $z_1$:
\begin{gather}
f(z_1,\dots,z_r)= f_0(z_2,\dots,z_r) +\sum_{i=1}^p \frac{\tilde f_i(z_2,\dots,z_r)}{y_i -z_1}.
\label{f}
\end{gather}
Then the function $f_0(z_2,\dots,z_r)$ has the properties $C_{r-p-1}(y_1,\dots,y_p;z_2,\dots,z_r)$, while
the function $\tilde f_i(z_2,\dots,z_r)$, ${i>0}$, has the properties $C_{r-p}(y_1,\dots,y_p;z_2,\dots,z_r)$ and
\[
\tilde f_i\big(q^{2} y_i,z_3,\dots,z_r\big) = 0,
\]
cf.~iii). The last claim is equivalent to the fact that the function
\begin{gather}
f_i(z_2,\dots,z_r) = \tilde f_i(z_2,\dots,z_r)\prod_{j=2}^r\frac{y_i -z_j}{q^{-1}y_i -q z_j}
\label{fi}
\end{gather}
has the properties $C_{r-p}(y_1,\dots,\widehat{y_i},\dots,y_p;z_2,\dots,z_r)$.

We expand the functions $f_0,\dots,f_p$ similarly to \eqref{f}, \eqref{fi}:
\[
f_i(z_2,\dots,z_r)= f_{i0}(z_3,\dots,z_r) +
\sum_{\substack{j=1 \\ j\ne i}}^p \frac{f_{ij}(z_3,\dots,z_r)}{y_j -z_2}
\prod_{s=3}^r \frac{q^{-1}y_i -q z_s}{y_i -z_s},
\]
and observe that the function $f_{00}$ has the properties $C_{r-p-2}(y_1,\dots,y_p;z_3,\dots,z_r)$,
the functions $f_{0i}$, $f_{i0}$, ${i>0}$ have the properties
$C_{r-p-1}(y_1,\dots,\widehat{y_i},\dots,y_p;z_3,\dots,z_r)$, and the function $f_{ij}$, ${i,j>0}$ has
the properties $C_{r-p}(y_1,\dots,\widehat{y_i},\dots,\widehat{y_j},\dots,y_p;z_3,\dots,z_r)$. Eventually, we obtain
the following expansion of the function $f(z_1,\dots,z_r)$:
\begin{gather}
f(z_1,\dots,z_r) = \sum_{\boldsymbol{\alpha}} f_{\boldsymbol{\alpha}}
\prod_{i=1}^r \left( \varphi_{\alpha_i}(z_i)
\prod_{\substack{i<j\leqslant r\\ \alpha_i>0}}
\frac{qy_{\alpha_i} -q^{-1} z_j}{y_{\alpha_i} -z_j} \right),
\label{pf}
\end{gather}
where the sum is taken over all surjective maps
$\boldsymbol{\alpha}:\{1,\dots,r\} \to\{0,\dots,p\}$ such that the preimage of $0$ has $r-p$
elements, $\varphi_0(u)=1$ and $\varphi_s(u)=(y_s -u)^{-1}$ for $i=1,\dots,p$.
The coef\/f\/icients $f_{\boldsymbol{\alpha}}$ do not depend on $z_1,\dots,z_r$ and can be found from
the equality
\begin{gather}
\operatorname{Val}_{\alpha_r,r} \cdots \operatorname{Val}_{\alpha_1,1}f(z_1,\dots,z_r) = (-1)^p q^{-c_{\boldsymbol{\alpha}}}
f_{\boldsymbol{\alpha}}\prod_{\substack{1\leqslant i<j\leqslant p \\ \alpha_i\alpha_j>0}}
\frac{q y_{\alpha_i} -q^{-1}y_{\alpha_j}}{y_{\alpha_i} -y_{\alpha_j}},
\label{Val}
\end{gather}
where ${\operatorname{Val}_{0,i}= \lim\limits_{z_i\to\infty}}$,
${\operatorname{Val}_{s,i}= \mathop{\operatorname{Res}}\limits_{z_i=y_s}}$ for ${s>0}$,  and
\[
c_{\boldsymbol{\alpha}}=
\#\big\{(i ,j)\,|\, {i<j}, \;{\alpha_i >0}, \;{\alpha_j =0} \big\}.
\]
Since the operations $\operatorname{Val}_{s,i}$ in the left side of \eqref{Val} can be applied to
the function $f(z_1,\dots,z_r)$ in any order without changing the answer,
$\operatorname{Val}_{\alpha_r,r}\cdots \operatorname{Val}_{\alpha_1,1}f(z_1,\dots,z_r)$ equals
\begin{gather}
\operatorname{Val}_{0,\tau_1}\cdots \operatorname{Val}_{0,\tau_p}
\operatorname{Val}_{1,\tau_{p+1}}\cdots \operatorname{Val}_{p,\tau_r}f(z_1,\dots,z_r)
\label{Valtau}
\end{gather}
for a suitable permutation $\tau$. Since $f(z_1,\dots,z_r)$ is symmetric in $z_1,\dots,z_r$,
expression~\eqref{Valtau} does not depend on $\tau$ and equals
\begin{gather}
\lim_{z_1\to\infty} \cdots\lim_{z_{r-p}\to\infty}
\mathop{\operatorname{Res}}\limits_{z_{r-p+1}=y_1} \cdots \mathop{\operatorname{Res}}\limits_{z_r=y_p}f(z_1,\dots,z_r).
\label{limres}
\end{gather}
Due to the explicit formula for $f(z_1,\dots,z_r)$, the terms in $\operatorname{Sym}_{z_1,\dots,z_r}$ which
contribute nontrivially to expression \eqref{limres} correspond to permutations that
do not move the numbers $r-p+1,\dots, r$. Using identity~\eqref{k!q}, we obtain
that expression~\eqref{limres} equals
\[
(-1)^p q^{-p(r-p)} [r-p]_q!\prod_{1\leqslant i<j\leqslant p}
\frac{\big(q^{-1}y_i -qy_j\big) \big(qy_i -q^{-1}y_j\big)}{(y_i -y_j)^2}.
\]
Hence, equality \eqref{Val} yields
\begin{gather}
f_{\boldsymbol{\alpha}} = q^{c_{\boldsymbol{\alpha}} -p(r-p)}[r-p]_q!
\prod_{\substack{1\leqslant i<j\leqslant p\\ \alpha_i\alpha_j> 0}}
\frac{q^{-1}y_{\alpha_i} -q y_{\alpha_j}}{y_{\alpha_i} -y_{\alpha_j}}.
\label{falph}
\end{gather}

There exists a bijection between pairs $({\boldsymbol{d}} ,\sigma)$, where ${\boldsymbol{d}}$ is a
$p$-tuple from Lemma~\ref{Gyzq} and $\sigma$ is a permutation of $\{1,\dots,p\}$,
and the maps $\boldsymbol{\alpha}$. It is given by the rule ${\alpha_{d_i}=\sigma_i}$, ${i=1,\dots,p}$,
and ${\alpha_j =0}$, otherwise. Under this bijection, the right side of formula~\eqref{pf}
with the coef\/f\/icients $f_{\boldsymbol{\alpha}}$ given by formula~\eqref{falph} turns into the right side
of formula~\eqref{yzq}.
\end{proof}

\begin{proof}[Proof of Lemma~\ref{Gyz}]
Make the change of variables  ${y_i\to 1+2h y_i}$, ${z_i\to 1+2hz_i}$,
$q\to 1+h$  in formula~\eqref{yzq} and take the limit $h \to 0$.
This yields the claim.
\end{proof}

\subsection*{Acknowledgments}

The authors thank the Max-Planck-Institut f\"ur Mathematik in Bonn
for hospitality.
V.~Tarasov was supported in part by NSF grant DMS--0901616.
A.~Varchenko was supported in part by NSF grants DMS--1101508.

\pdfbookmark[1]{References}{ref}
\LastPageEnding

\end{document}